\documentclass[leqno]{amsart}

\usepackage[sort&compress]{natbib}

\usepackage[foot]{amsaddr}

\usepackage{color}
\usepackage{geometry}
\usepackage{subfiles}
\usepackage{graphicx}

\usepackage{caption}
\usepackage{subcaption}
\usepackage{float}

\usepackage[section]{placeins}
\usepackage{hyperref}
\usepackage{stmaryrd}
\usepackage{amsmath,amssymb,amsfonts,amsthm,textcomp}
\usepackage{mathtools}
\usepackage{tensor}
\usepackage{multicol}
\usepackage{enumitem}
\usepackage{tcolorbox}
\usepackage{bbold}
\usepackage{tikz}
\usepackage{array}
\usepackage{longtable}

\usetikzlibrary{spath3, intersections, hobby, knots}
\usepackage{pst-knot}
\usepackage{spath3}
\usepackage{enumitem}

\graphicspath{
              {./Figures/}
             }

\newfont{\tenbbb}{msbm10}
\newfont{\svnbbb}{msbm8}

\usepackage[T3,T1]{fontenc}
\DeclareSymbolFont{tipa}{T3}{cmr}{m}{n}
\DeclareMathAccent{\invbreve}{\mathalpha}{tipa}{16}

\theoremstyle{remark}

\theoremstyle{definition}

\makeatletter
\def\threevdots{\mskip+4mu\vbox{\baselineskip2.25\p@ \lineskiplimit\z@
  \kern4.9\p@\hbox{.}\hbox{.}\hbox{.}}\mskip+3.8mu}
\makeatother

\theoremstyle{definition}
\newtheorem{definition}{Definition}[section]
\newtheorem{theorem}{Theorem}[section]
\newcommand{\knotlinewidth}{.7pt}
\newcommand{\myarrow}{$\quad\longleftrightarrow\quad$}
\newcommand{\RIa}[1][]{\tikz[knot, #1]{\draw(-.5,.5) to[out=-90,in=-90] (.5,0); \draw[overcross] (.5,0) to[out=90,in=90] (-.5,-.5);}}
\newcommand{\RIb}[1][]{\tikz[knot, #1]{\draw[looseness=.8] (-.5,-.5) to[out=90, in=-90] (.5,0) to[out=90, in=-90] (-.5,.5);}}
\newcommand{\RIIa}[1][]{\tikz[knot, #1]{\draw[red, looseness=2.3] (-.5,-.5) to[out=0, in=0] (-.5,.5); \draw[looseness=2.3, overcross=blue] (.5,-.5) to[out=180, in=180] (.5,.5);}}
\newcommand{\RIIb}[1][]{\tikz[knot, #1]{\draw[red, looseness=1.4] (-.5,-.5) to[out=0, in=0] (-.5,.5);\draw[blue, looseness=1.4] (.5,.5) to[out=180,in=180] (.5,-.5);}}
\newcommand{\RIIc}[1][]{\tikz[knot, #1]{\draw[blue, looseness=2.3] (.5,-.5) to[out=180, in=180] (.5,.5); \draw[looseness=2.3, overcross=red] (-.5,-.5) to[out=0,in=0] (-.5,.5);}}
\newcommand{\RIIIa}[1][]{\tikz[knot, #1]{\draw[red] (-120:.58) to[out=60,in=-120] (150:.2) to[out=60, in=-120] (60:.58); \draw[rotate=-120, overcross=blue] (-120:.58) to[out=60, in=-120] (150:.2) to[out=60, in=-120] (60:.58); \draw[rotate=120, overcross=green!80!black] (-120:.58) to[out=60, in=-120] (150:.2) to[out=60, in=-120] (60:.58);}}

\tikzset{overcross/.style={double, line width=1.5, white, double=#1, double distance=\knotlinewidth},
    overcross/.default={black},
    knot/.style={line width=\knotlinewidth, baseline=-.5ex}}

\begin{document}

\title[All Tied Up: A first look at modern Yo-Yo and knot theory]{All Tied Up: A first look at modern Yo-Yo and knot theory}
\author{Benjamin Hamblin$^\sharp$, Victor M. Calo$^\sharp$}

\address{$^\sharp$School of Electrical Engineering, Computing and Mathematical Sciences, Curtin University, Bentley, Western Australia}
\email{benjamin.hamblin@curtin.edu.au}

\date{\today}

\begin{abstract}
\noindent
Modern Yo-Yo play has developed into a sophisticated international subculture, featuring elite competition and intricate tricks. Despite this, no systematic knot-theoretic treatment has yet been applied to the numerous string configurations realised in contemporary Yo-Yo play. This paper takes initial steps in addressing this gap, recalling fundamental results from knot theory and embeddings to develop a methodology for classifying string arrangements, known as `mounts', in both beginner and advanced single Yo-Yo play.  
We classify a range of mounts according to the knots they form under appropriate post-processing procedures and identify Yo-Yo maneuvers that correspond to Reidemeister moves. Furthermore, we analyse the impact of certain mounts on the writhe of their diagrammatic projections and introduce operations to facilitate discussion of composite mounts.  
This work seeks to initiate a dialogue between Yo-Yo practitioners and knot theorists, fostering further advancements in high-level Yo-Yo play and enabling novel physical realisations of knots, links, braids, surgeries, and other topological transformations.

\end{abstract}

\maketitle
\tableofcontents                      


\section{Introduction and Motivation}
The study of puzzles, hobbies, and crafts via mathematical framing is a rich tradition and serves both to further mathematical techniques and provide innovative avenues for education. For example, the intriguing geometry of origami has received considerable attention in recent times and is useful in applications such as robotics, biomedicine, and industrial design~\cite{ Origami}. Sudoku is used as an algorithmic testing ground and has strong connections with group theory~\cite{ Sudoku}. The combinatorics and notoriously P-hard strategy of Go made international headlines in 2015 when DeepMind's AI model AlphaGo defeated a top-ranked human player for the first time~\cite{ Go}. Perhaps surprisingly, the Yo-Yo has received no topological or knot-theoretic treatment to the author's knowledge, despite Yo-Yo play being a highly developed subculture with a large international following. There exist many investigations into the physics of Yo-Yo mechanics (e.g., see,~\cite{Yoyophysics1, minkin2020yo, yoyophysics2}). Yet, any formal analysis of the intricate string tricks that modern Yo-Yo players play is lacking. The purpose of this work is to open the door to such an analysis and develop a system of categorising string configurations in Yo-Yo play built on concepts borrowed from knot theory. This approach is hoped to be useful or interesting to at least three groups of people: 
\begin{enumerate}
\item \textbf{Yo-Yo players}, since formally categorising the crossings and twists in a given string arrangement (a ``mount'') offers a systematic framework to identify feasible transitions and to perform complex tricks without creating knots.  
\item \textbf{Mathematicians}, as the Yo-Yo provides a tangible means of realising intricate topological arrangements and exploring forms of ``knot surgery.''  
\item \textbf{Mathematics educators}, who can use a popular toy to introduce key concepts in geometry and topology, offering an engaging entry point to topics that are often perceived as abstract or difficult.  
\end{enumerate}
The remainder of this work is structured as follows. In section~\ref{History}, we discuss the highly developed phenomenon of modern competitive Yo-Yo play and introduce some basic terminology from that community relevant to the descriptions provided herein. In section~\ref{Knot Theory}, we give a summary of important notions from knot theory that we utilise in section~\ref{1A} to develop a systematic way of classifying Yo-Yo mounts (subject to certain constraints) according to their knot-theoretic properties, and apply this method to a range of conventional and intricate mounts. We conclude and outline potential avenues for further work in the area in section~\ref{Conclusion}.
\section{Modern Yo-Yo: Fundamentals and History}
\label{History}
The Yo-Yo is an ancient toy, with origins traceable to at least 440\,BC. In recent decades, advances in design have transformed it from a simple pastime into a platform for remarkably complex tricks across multiple styles of play. This evolution has given rise to a thriving global subculture of competitive Yo-Yo players. The highest stage of this culture is the annual World Yo-Yo Contest (WYYC), overseen by the International Yo-Yo Federation (IYYF), which became independent of the International Juggling Association in 1998.  

Technological innovations, such as widened gaps and concave bearings, have enabled the development of intricate string tricks using modern \emph{unresponsive} Yo-Yos—designs that do not return with a simple tug, but instead require a specialised trick, known as a \emph{bind}, to return the Yo-Yo to the hand. In contemporary competitions, players perform three-minute freestyles (in final rounds), choreographing elaborate sequences of tricks to music. Performances are judged on two criteria~\cite{iyyf}. The first is \emph{Technical Execution} (TE), in which points are awarded for each non-trivial element, with additional points for more difficult maneuvers. The second is \emph{Freestyle Evaluation} (FE), in which judges assess performance qualities such as control, creativity, showmanship, and overall presentation. \\ 

\begin{center}
\begin{figure}[h!]
\caption{A modern 1A Yo-Yo mount}
\includegraphics[scale=0.25]{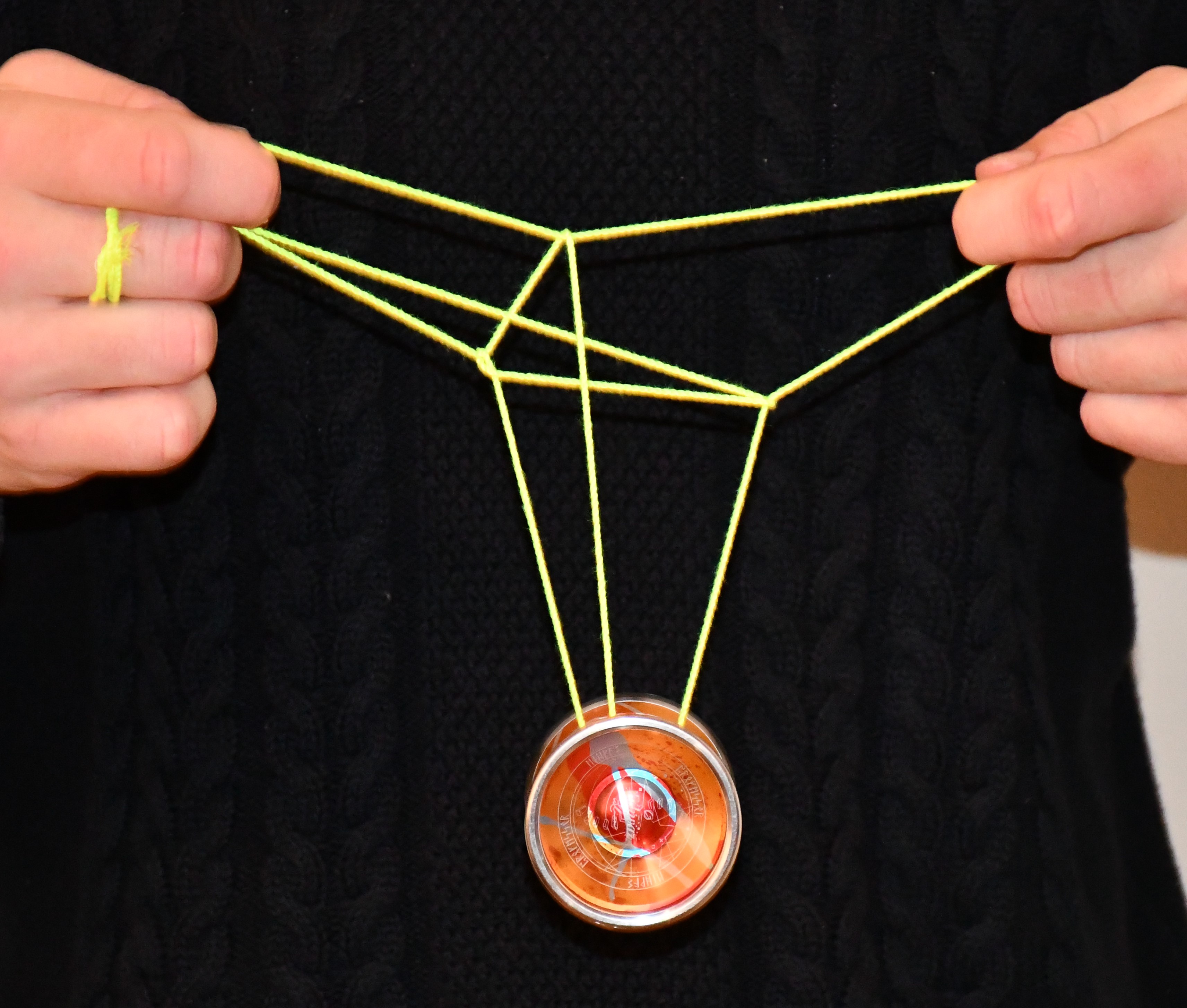}
\label{Complex}
\end{figure}
\end{center}

For elite players, freestyles are taken seriously, as prestige and lucrative sponsorships are at stake. While many in the community share tutorials for developing players, the creative process behind new tricks is often treated as tradecraft, with signature routines typically revealed only on the contest stage. Despite this, most tricks share a common structure: transitions between configurations of string known as \emph{mounts}. A mount arranges the string so that the Yo-Yo rests on one or more segments. There exist a vast number of possible mounts, many of which form intricate and visually striking structures (see Figure~\ref{Complex}).  

Certain mounts, however, can appear highly complex yet reduce to a much simpler configuration with a single move. Knot theory provides a natural framework to characterise this complexity, measuring it in terms of the number of non-trivial crossings. In this way, it offers insight into why intricate mounts can collapse into simpler ones—a phenomenon we aim to demonstrate in the remainder of this work.  

\subsection{Styles of Modern Competition Yo-Yo}
Modern Yo-Yo play encompasses a wide range of distinct styles, with new variations continually emerging as players push the boundaries of creativity. To illustrate both the breadth and sophistication that contemporary Yo-Yo has achieved, we briefly outline the five most common styles, which also form the official championship divisions at the World Yo-Yo Contest: 

\begin{enumerate}
\item \textbf{1A: Single Yo-Yo String Tricks}\\
Currently the most popular competitive style, 1A involves manipulating a single unresponsive Yo-Yo and string into intricate mounts and transitions. Typical elements include slacks, whips, hops, and rejections, many of which form the foundation for highly elaborate routines. 

\item \textbf{2A: Two-Handed Looping Tricks}\\
This fast-paced style uses two responsive Yo-Yos simultaneously. Players perform continuous looping patterns with both hands, creating a dynamic and rhythmic display that demands speed, precision, and endurance.

\item \textbf{3A: Double Yo-Yo String Tricks}\\
Regarded as one of the more technically demanding divisions, 3A features two unresponsive Yo-Yos manipulated simultaneously. The interaction between the two Yo-Yos enables the creation of highly intricate string structures and maneuvers beyond the scope of 1A play.  

\item \textbf{4A: Offstring Yo-Yo}\\
In 4A, the Yo-Yo is not tied to the string but instead launched and caught repeatedly, enabling aerial tricks such as orbits and whips. The style borrows from Diabolo and includes advanced sub-styles such as Soloham, in which players manipulate two off-string Yo-Yos at once.  

\item \textbf{5A: Freehand Yo-Yo}\\
Here, the string is not anchored to the player’s finger but attached to a counterweight, typically about one-sixth the mass of the Yo-Yo. This frees both hands for manipulation, making possible dynamic tricks that draw inspiration from other performance arts such as Poi and Meteor.  
\end{enumerate}

For conciseness, our subsequent analysis will focus primarily on the 1A style. Nevertheless, many of the principles we discuss extend naturally to mounts and transitions found in other divisions.  

\section{Knot Theory: Some Useful Concepts}
\label{Knot Theory}
\subsection{Equivalent Knots and Reidemeister Moves}
A knot is an embedding of the circle \(S^1\) into Euclidean three-space, which we denote \(\mathbb{E}^3\). 
\begin{definition}[Knot Equivalence]
Two knots \(k_1, k_2\) are said to be equivalent up to ambient isotopy if there exists an orientation-preserving homeomorphism \(h\) of \(\mathbb{E}^3\) such that \(h(k_1) = k_2\).
\end{definition}
A crossing is any point on a diagram of a knot where the knot intersects itself. 
\begin{definition}[Crossing Number]
The crossing number of a knot \(k\) is the smallest number of crossings possible on a diagram of a knot equivalent to \(k\).
\end{definition}
The crossing number provides a useful first measure of knot complexity, but on its own it cannot fully capture when two knots are equivalent. To build a more powerful framework, we turn to the Reidemeister moves. These three fundamental moves encompass all possible local transformations of a knot diagram. When combined with isotopies (continuous deformations such as translations or smooth bending of the string without cutting or gluing), they form a complete toolkit: any two diagrams of the same knot can be transformed into one another by a finite sequence of Reidemeister moves. \\
The moves are:
\begin{enumerate}[label=(\Roman*)]
\item Twist and untwist in either direction\\
\RIa\myarrow\RIb\myarrow\RIa[yscale=-1]
\item Move one loop completely over another\\
\RIIa\myarrow\RIIb\myarrow\RIIc
\item Move a string completely over or under a crossing\\
\RIIIa\myarrow\RIIIa[rotate=180]\qquad\RIIIa[xscale=-1]\myarrow\RIIIa[xscale=-1, rotate=180]
\end{enumerate}

We can use Reidemeister moves to invoke a more practical notion of knot equivalence, specifically:
\begin{theorem}[Reidemeister's Theorem]
Two knots are equivalent if and only if they result from one another via a finite series of Reidemeister moves. 
\end{theorem}
\textit{Proof}\hspace{0.5cm} The proof follows as per Reidemeister's original proof in~\cite{Reidemeister}. \qed\vspace{0.1cm}\\

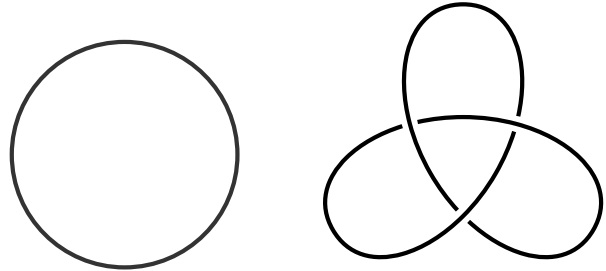
\begin{figure}[h!]

\caption{Diagrams of the unknot and the simplest non-trivial knot, the trefoil $(3_1)$}

\begin{tikzpicture}[
  use Hobby shortcut,
  every trefoil component/.style={ultra thick, draw, black},
  trefoil component 1/.style={black},
]

\draw[color=black!80, ultra thick](-4,0) circle (1.5);
\hspace{0.5cm}
\path[spath/save=trefoil] ([closed]90:2) foreach \k in {1,...,3} { .. (-30+\k*240:.5) .. (90+\k*240:2) } (90:2);
\tikzset{spath/knot={trefoil}{6pt}{1,3,5}}

\end{tikzpicture}
\label{knots}
\end{figure}

This is an important result for the topic herein, as it demonstrates that two arrangements of a closed loop of string may be in some sense the same arrangement despite appearing very different. With this toolkit in place, mathematicians have created a `dictionary of knots' wherein the distinct (non-equivalent) prime knots of a given crossing number are indexed according to a subscript. A prime knot is a knot that cannot be written as the knot sum of two non-trivial knots, a notion we discuss in more detail in~\ref{Knot Sum}. Since there are infinitely many prime knots~\cite{KnotBook}, naturally, this `dictionary' is not exhaustive, but there exist comprehensive listings for sufficiently small crossing numbers. The unknot is the trivial knot with no crossings, and will be denoted for the duration of this work by \(\mathfrak{u}\). The simplest non-trivial knot is the trefoil, which contains three crossings, and is denoted by \(3_1\). Diagrams of these knots may be seen in Figure~\ref{knots}. 

\subsection{The Knot Sum and Prime Knots}
\label{Knot Sum}
\begin{definition}[Connected Knot Sum]
The \textit{connected knot sum}, or \textit{composition}, of two knots \(k_1\) and \(k_2\) is the knot obtained by the attachment of \(k_1\) to \(k_2\) with respect to the orientation of both knots.
\end{definition}
Hereafter, we call the connected knot sum simply the `knot sum' for brevity. The result of `attaching' \(k_1\) to \(k_2\) is obtained via the following procedure~\cite{KnotBook2}:
\begin{enumerate}
\item Locate an arc that includes no crossing in \(k_1\) and an arc that includes no crossing in \(k_2\) that are adjacent to each other in the sense that there exist continuous paths from the endpoints of one to the other without intersecting any part of either knot. 
\item Delete these arcs from their corresponding knots. Call what is left of \(k_1\), \(k_1'\) and call what is left of \(k_2\), \(k_2'\).
\item Connect the end-points of \(k_1'\) with \(k_2'\) in such a way that no additional crossings are introduced. 
\end{enumerate}
We call the knot arising from this procedure \(k_1\#k_2\). 
\begin{definition}[Prime Knot]
A \textit{prime knot} is a non-trivial knot which cannot be written as the knot sum of two non-trivial knots. 
\end{definition}
Oriented knots in three-dimensional space form a commutative monoid under the knot sum, with the unknot serving as the identity element, and they admit a unique prime factorisation~\cite{Manifolds}. Orientation plays a central role in our setting, since the direction in which the Yo-Yo travels (under appropriate procedures for recovering mathematical knots) may yield, for instance, either a left-handed or right-handed trefoil. These two are not equivalent: no finite sequence of Reidemeister moves can transform one into the other~\cite{KnotBook}. In this sense, the trefoil is said to be \textit{chiral}, meaning that it is distinct from its mirror image.  

Topological orientation is a subtle matter, and several notational conventions have been proposed to capture it. In what follows, we adopt a straightforward scheme: we denote left-handed and right-handed orientations with subscripts \(L\) and \(R\), respectively. Later, in Section~\ref{Setup}, we introduce a practical notion of orientation as it arises in the context of Yo-Yo mounts and, in special cases, compare it explicitly to handedness. For the moment, we collect a few key observations:
\begin{equation}
3_{1,L} \# 3_{1,R} = \mathfrak{u}, \hspace{0.5cm}\textrm{but}\hspace{0.5cm} 3_{1,X} \# 3_{1,X} \ne \mathfrak{u},
\end{equation}
where \(X\) may denote \(L\) or \(R\). That is to say, the knot sum of a trefoil with its mirror image is the unknot, whilst the sum of a given oriented trefoil with a copy of itself is a composite knot with six crossings, known as the `granny knot'. 
\subsection{Writhe}
A useful notion to distinguish between equivalent knots\footnote{Distinguishing between, say, a trefoil and a trefoil with an additional twist may be essential depending on the context. For instance, certain DNA helices can be topologically equivalent to the unlink, yet the extent to which they are \emph{coiled} plays a fundamental role in describing their mechanical and biological behaviour~\cite{DNA}.}, whether represented as planar diagrams or as closed curves embedded in \(\mathbb{E}^3\), is the so-called \textit{writhe}. For embedded curves, the writhe takes the Gauss integral expression~\cite{Writhe},
\begin{equation}
\label{writhe}
\textrm{Wr} = \frac{1}{4\pi}\int_C\int_C \textrm{d}\mathbf{r}_1 \times \textrm{d}\mathbf{r}_2\cdot\frac{\mathbf{r}_1 - \mathbf{r}_2}{|\mathbf{r}_1 - \mathbf{r}_2|^3},
\end{equation}
where \(C\) is the embedded curve, and \(\mathbf{r}_1, \mathbf{r}_2\) are points along \(C\). 
Intuitively, the writhe measures the total signed self-entanglement of the curve by accounting for all pairwise crossings of infinitesimal curve elements, weighted by their geometric orientation. In diagrammatic terms, it reduces to the signed sum of crossings in a knot projection, but the integral formulation emphasises its geometric rather than purely combinatorial nature. This expression provides a measure of how `twisted' the knot is. The integral in~\eqref{writhe} is equivalent to the integral average of all possible \textit{diagrammatic} writhes of a given curve, which takes integer values for a given projected knot (or link) diagram. For simplicity, it is the latter which we will concern ourselves with, although it is worth mentioning that mounts which are non-planar in construction may exhibit in-play mechanics which are more appropriately described by the embedded curves notion of writhe. 
\begin{definition}[Diagrammatic Writhe]
The \textit{writhe} (or diagrammatic writhe) for a given knot diagram is the total number of positive crossings minus the total number of negative crossings.
\end{definition}
We define the orientation of the strings used to construct knots in Section~\ref{1A} to run from the origin toward the terminus of the string. With this convention, a \textit{positive crossing} is one in which the overhand strand runs from left to right, whereas a \textit{negative crossing} occurs when the overhand strand runs from right to left.  

Within this framework, the \emph{writhe} emerges as a useful descriptor. Although two mounts may be topologically equivalent---in the sense that one can be transformed into the other by a sequence of Reidemeister moves---their writhe can differ. This difference encodes information about the way the mounts are constructed, which in turn may influence how they behave during play.
\section{1A Yo-Yo and Knot Theory}
\label{1A}
\subsection{Setup and Assumptions}
\label{Setup}

A Yo-Yo mount is not a physical realisation of a mathematical knot, because the string terminates at the Yo-Yo bearing rather than forming a closed loop with the origin at the player's finger. We address this by constructing two knots from a given mount via the following procedure:

\begin{enumerate}
    \item Disregard the slipknot placed over the player's finger, and imagine that the origin of the string is an ordinary segment.\footnote{There exists an unconventional style of Yo-Yo play, called `Moebius', in which the slipknot is topologically significant. In this style, the player may expand the slipknot and pass the Yo-Yo through it, making the slipknot integral to the mount. Here, we restrict attention to conventional 1A mounts.}
    \item Once the Yo-Yo is placed in a physical mount, imagine removing it and suspending the string in formation in mid-air.
    \item Construct the first knot by moving the terminus of the string \emph{in front of} the segment it is currently in contact with (from the player's point of view) and connect it to the origin. 
    \item Construct the second knot by moving the terminus \emph{behind} the segment it is currently in contact with and connect it to the origin. 
\end{enumerate}

The two procedures of merging the terminus with the origin—corresponding to the first and second knots—are referred to as \textit{hopping forward} and \textit{hopping backward}. These produce physical realisations of mathematical knots, giving a forward-hopped and a backward-hopped knot for each mount.  
Ambiguities may arise when the Yo-Yo bearing contacts more than one string segment besides the terminus. In such cases, the forward-hopped knot is obtained by moving the terminus \emph{in front of all} contacted segments, and the backward-hopped knot \emph{behind all} segments. An example is the `Kamikaze' mount (see Figure~\ref{Kamikaze}).  

\begin{definition}[Mount]
A mount is a double of knots.
\end{definition}

We denote forward-hopped knots with a \(\sharp\) superscript and backward-hopped knots with a \(\flat\) superscript. For example:
\begin{equation}
\label{Mount Examples}
T = (\mathfrak{u}^\sharp, \mathfrak{u}^\flat), \hspace{1cm} GT^\sharp = (\mathfrak{u}^\sharp, 3_1^\flat),
\end{equation}
where the superscripts are technically redundant, as the forward-hopped knot is listed first by convention, but are included for clarity. These correspond to the Trapeze and the Green Triangle (Forward-Hopped).\footnote{We include `Forward-Hopped' in brackets because most Yo-Yo players would call this mount simply a Green Triangle. For our purposes, there is a distinction between a Green Triangle yielding the unknot when hopped forward versus backward.}

In~\eqref{Mount Examples}, \(T\) forms the unknot when hopped forward and backward, while \(GT^\sharp\) forms the unknot when hopped forward and the \(3_1\) (trefoil) knot when hopped backward. Here, `forms the knot' means that the knot arises when disregarding the slipknot, removing the Yo-Yo, and connecting the terminus to the origin without introducing additional crossings.

\subsection{1A Mounts Admitting Knot-Theoretic Descriptions}
\label{Knot Mounts}

\begin{figure}[h!]
    \centering
    \begin{subfigure}{\textwidth}
        \centering
        \includegraphics[width=0.5\textwidth]{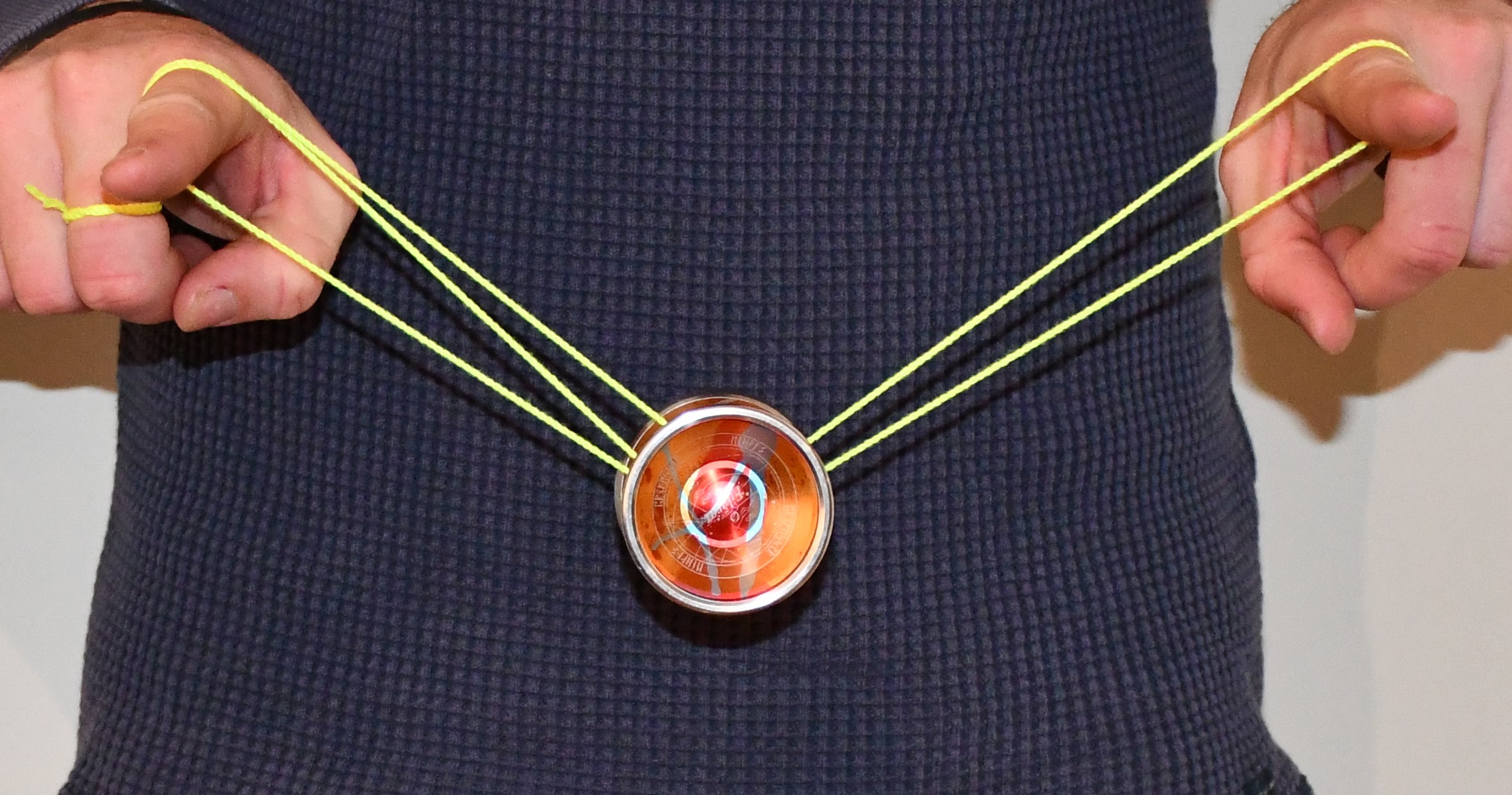}
        \caption{The Kamikaze mount in play}
        \label{fig:Kamikaze_in_play}
    \end{subfigure}\\
    \vspace{.2cm} 
    \begin{subfigure}[b]{0.35\textwidth}
        \centering
        \includegraphics[width=\textwidth]{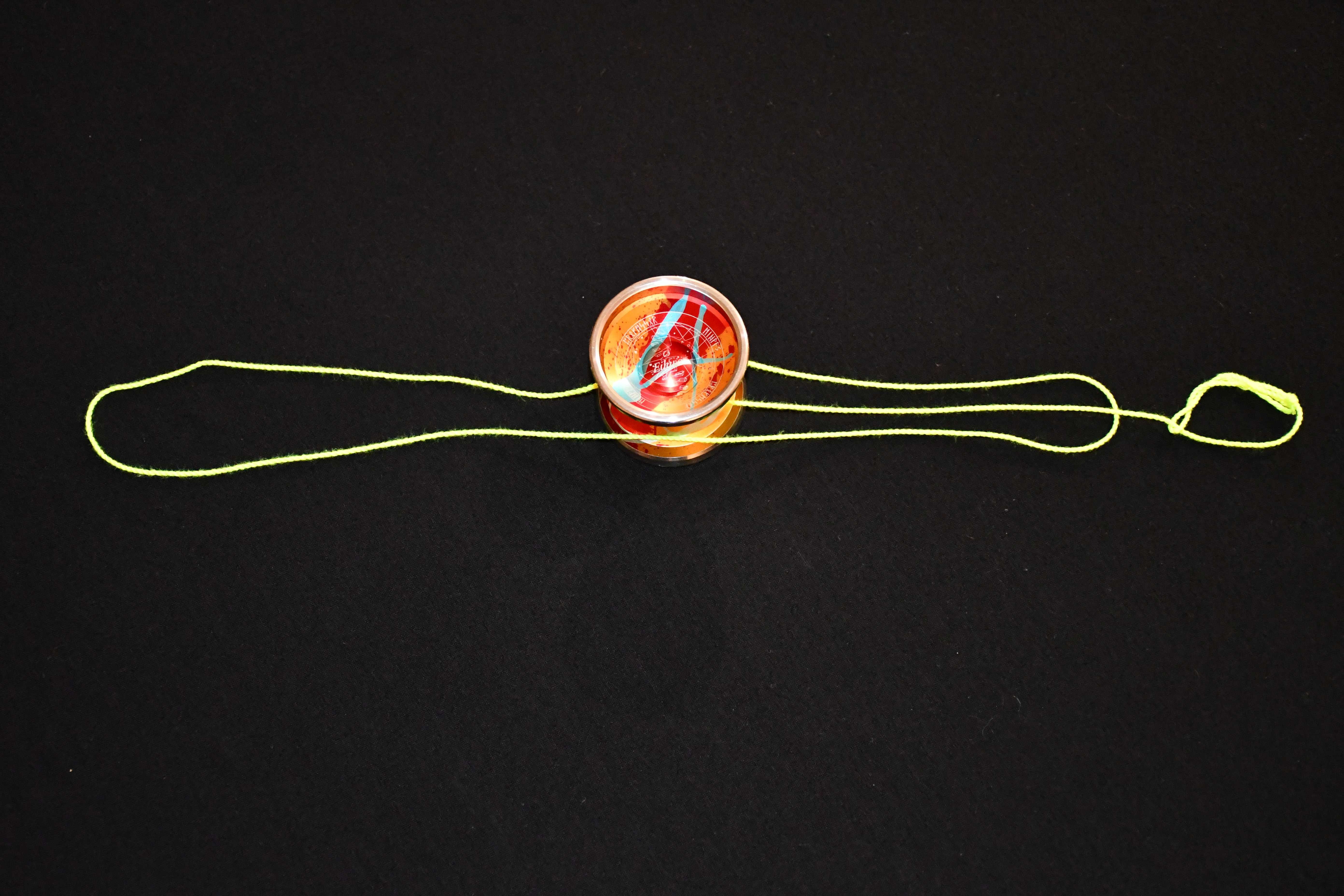}
        \caption{Mount with a stationary Yo-Yo}
        \label{fig:Kamikaze_flat}
    \end{subfigure}
    \quad
    \begin{subfigure}[b]{0.35\textwidth}
        \centering
        \includegraphics[width=\textwidth]{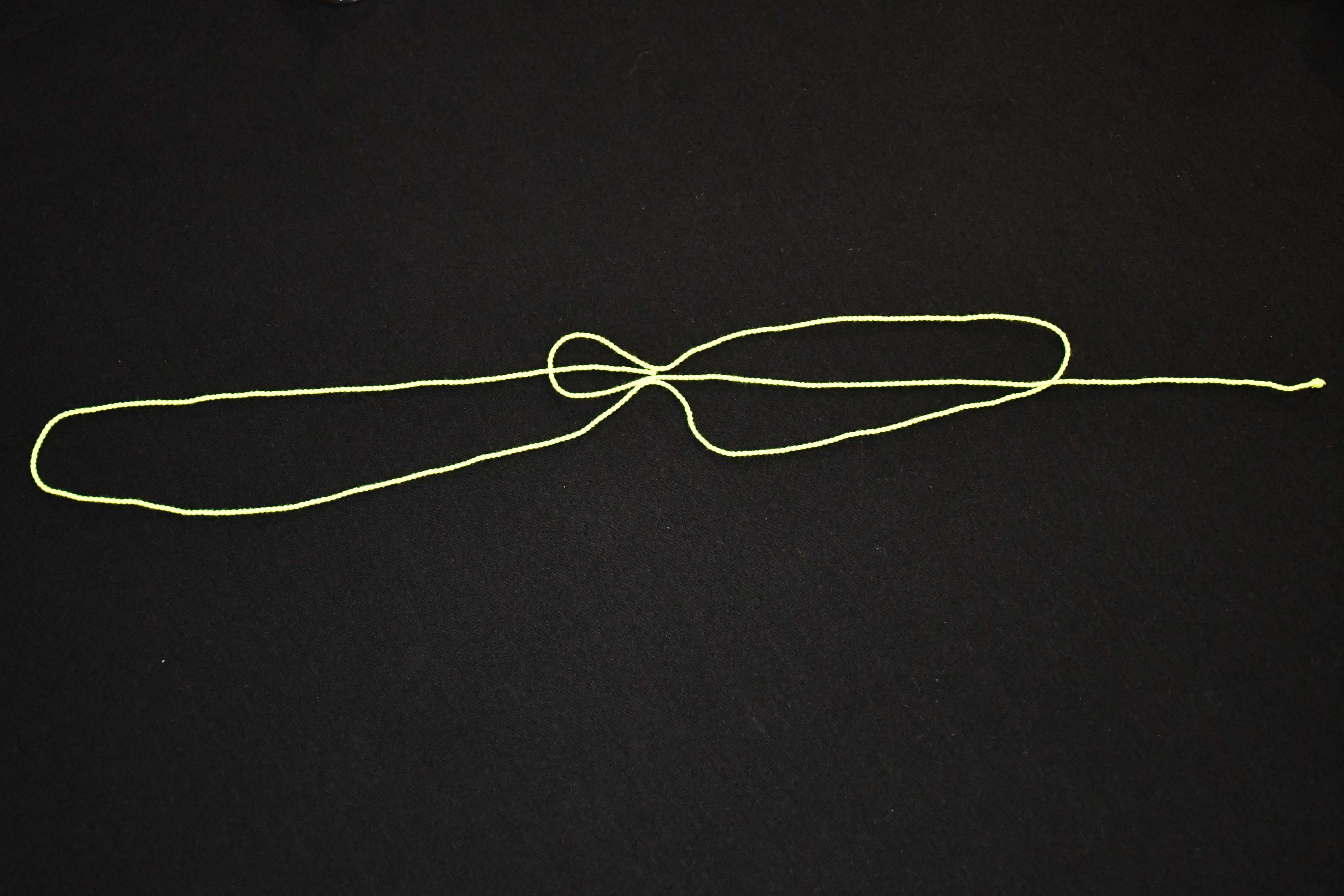}
        \caption{Slipknot replaced with a string segment}
        \label{fig:Kamikaze_no_yoyo}
    \end{subfigure}
    \\
    \vspace{.1cm} 
    \begin{subfigure}[b]{0.35\textwidth}
        \centering
        \includegraphics[width=\textwidth]{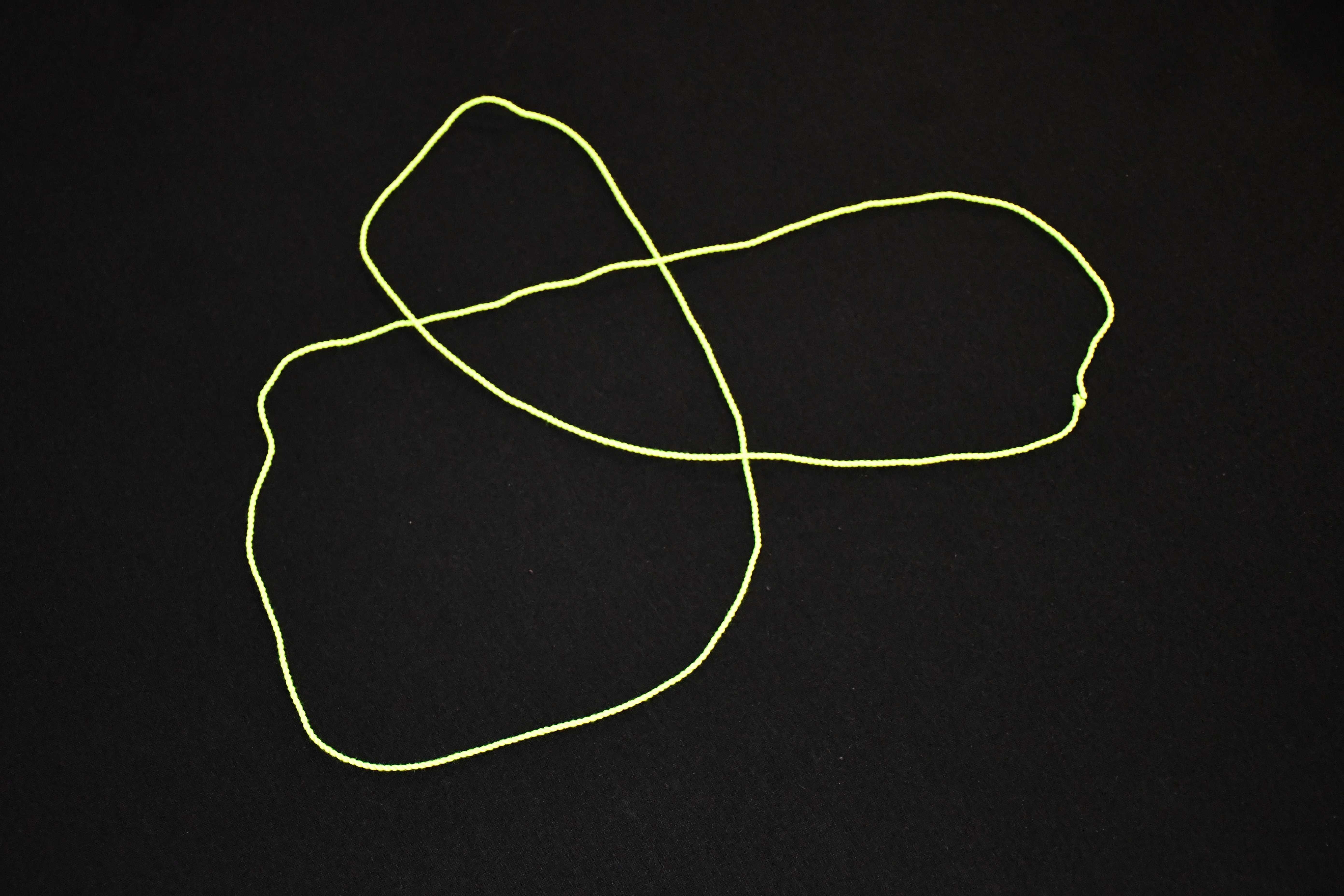}
        \caption{Forward-hopped knot}
        \label{fig:Kamikaze_forward}
     \end{subfigure}
     \quad 
    \begin{subfigure}[b]{0.35\textwidth}
        \centering
        \includegraphics[width=\textwidth]{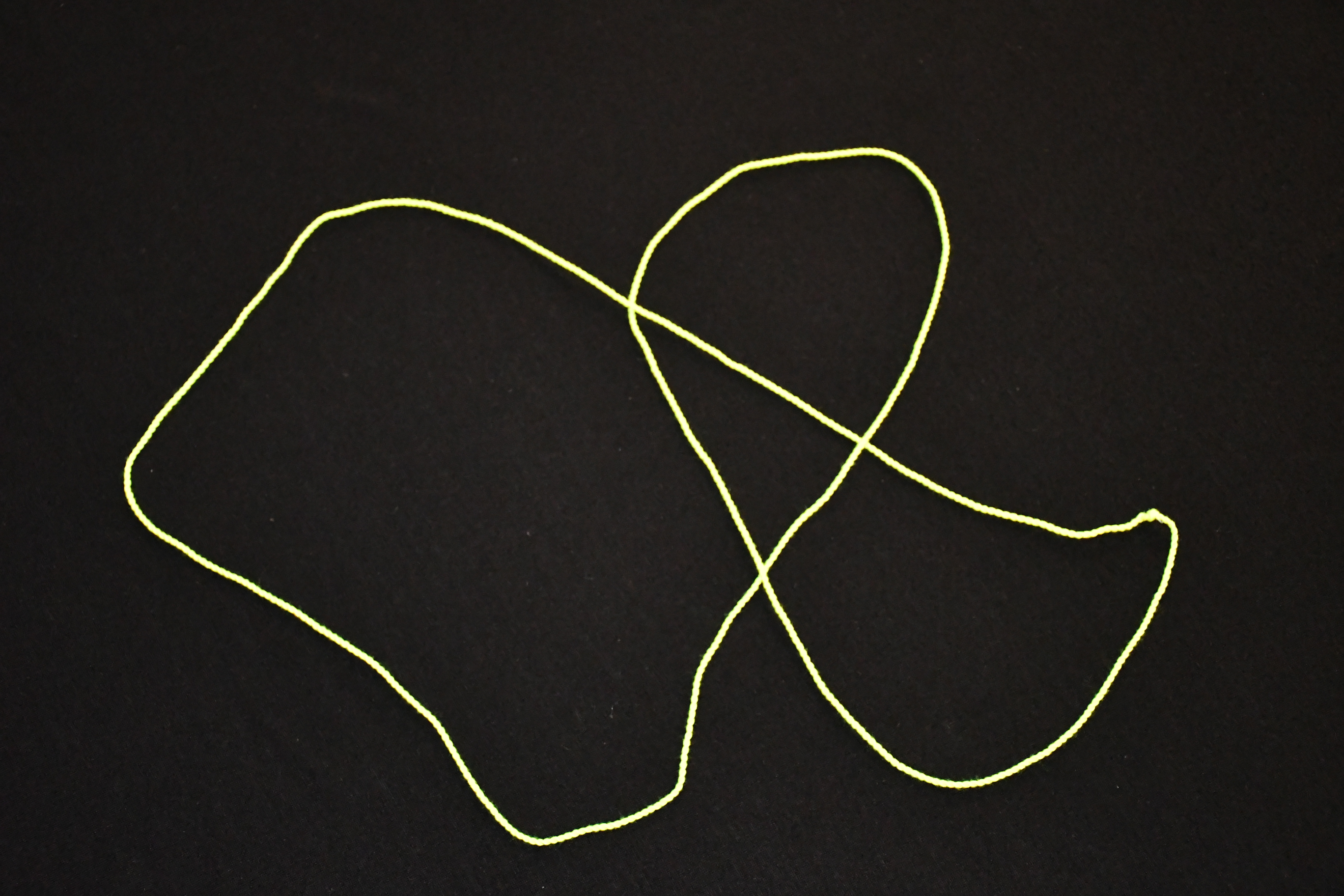}
        \caption{Backward-hopped knot}
        \label{fig:Kamikaze_backward}
    \end{subfigure}
    \caption{
      Construction of the forward- and backward-hopped knots for the Kamikaze mount. 
      \textbf{Top:} The mount in play. 
      \textbf{Bottom Rows:} Key stages of the procedure: the mount lying flat with the Yo-Yo stationary (\subref{fig:Kamikaze_flat}), the mount with the slipknot replaced by a simple string segment (\subref{fig:Kamikaze_no_yoyo}), and the resulting forward-hopped (\subref{fig:Kamikaze_forward}) and backward-hopped (\subref{fig:Kamikaze_backward}) knots. Both final knots recover the unknot within two Reidemeister moves.
    }
    \label{Kamikaze}
\end{figure}

\begin{figure}[H]
    \centering
    \begin{subfigure}[b]{0.45\textwidth}
        \centering
        \includegraphics[width=0.5\linewidth]{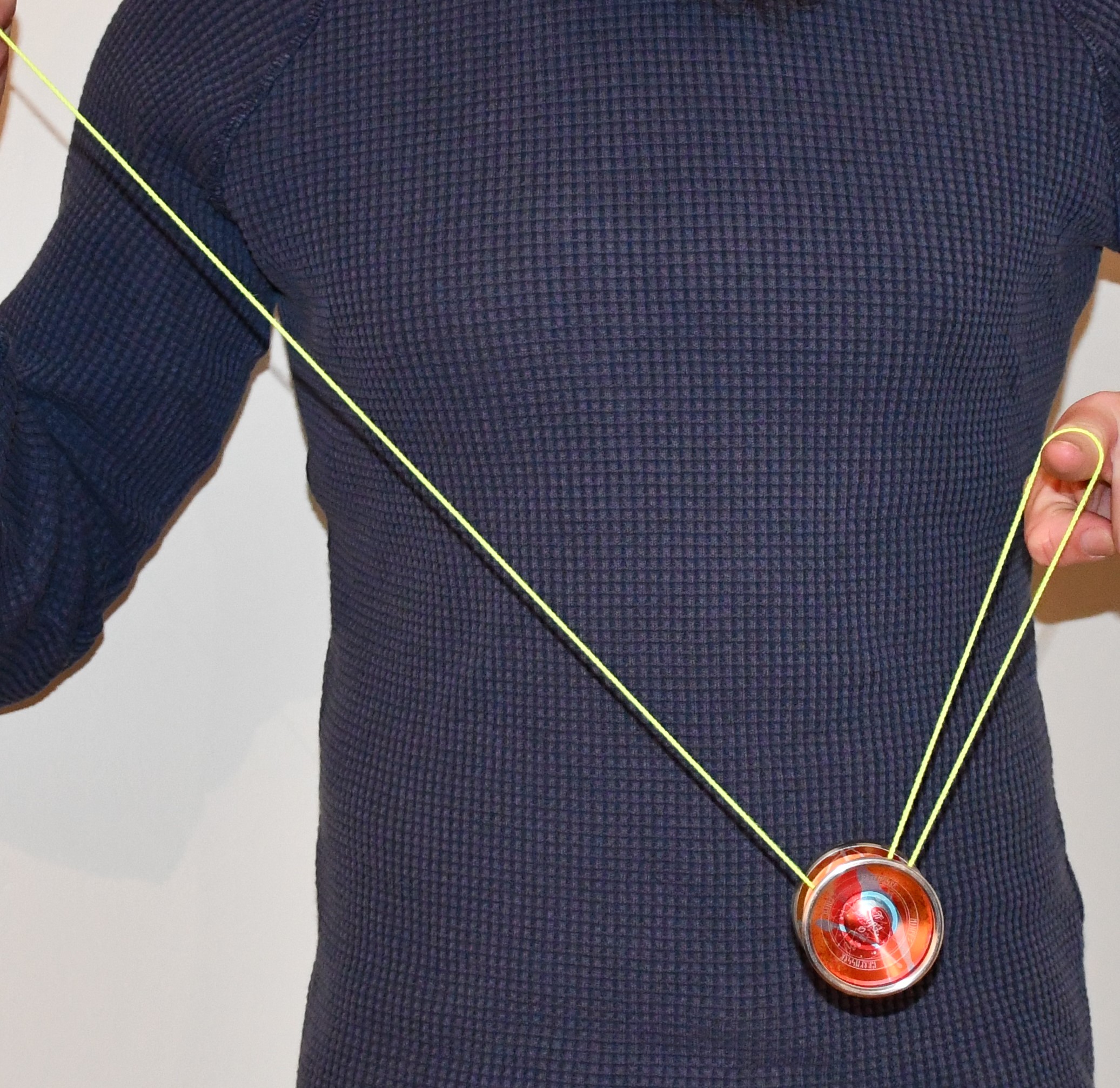}
        \caption{Trapeze}
        \caption*{\(T = (\mathfrak{u}^\sharp, \mathfrak{u}^\flat)\)}
        \label{fig:trapeze}
    \end{subfigure}
    \qquad
    \begin{subfigure}[b]{0.45\textwidth}
        \centering
        \includegraphics[width=0.6\linewidth]{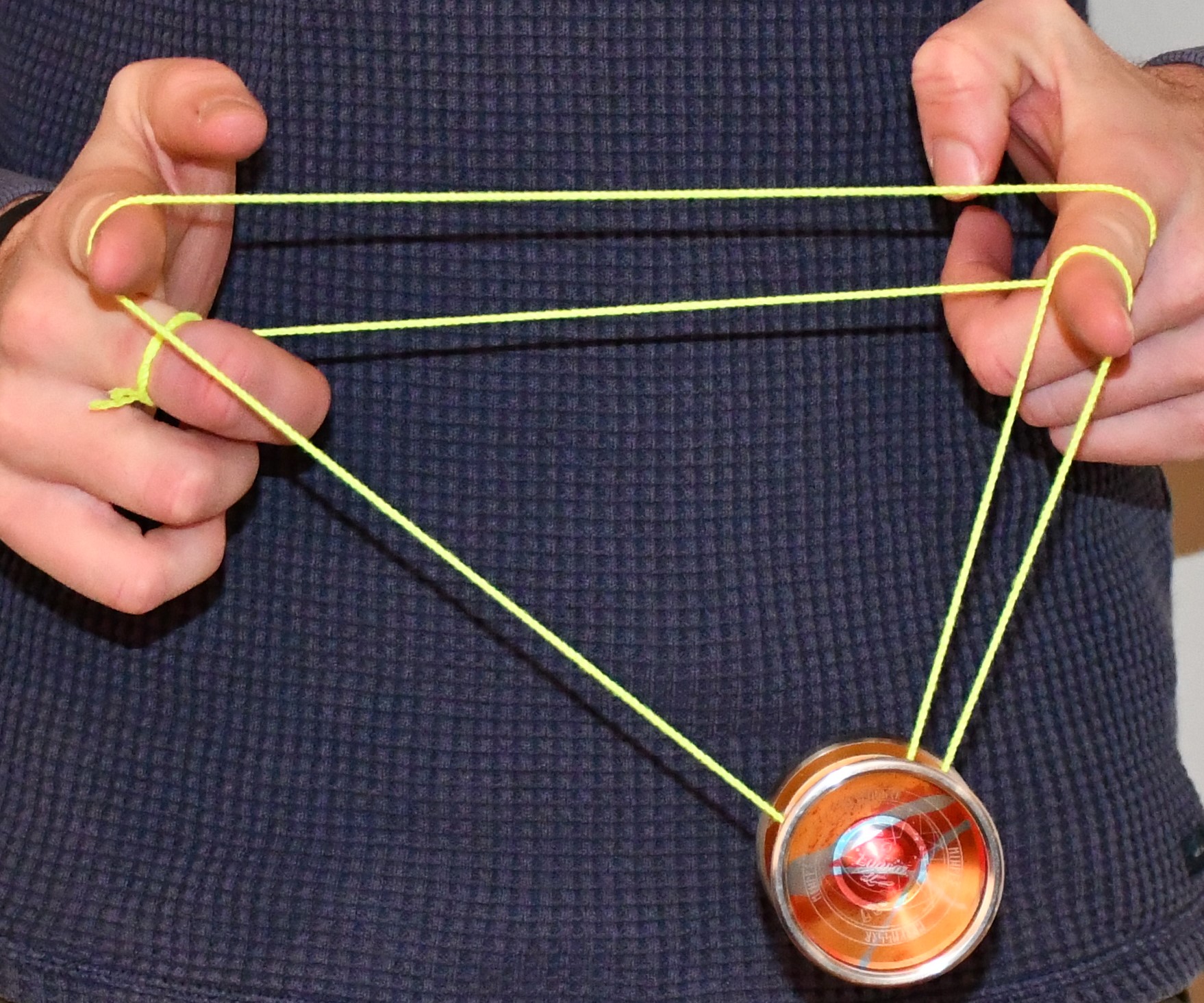}
        \caption{Double or Nothing}
        \caption*{\(D = (\mathfrak{u}^\sharp, \mathfrak{u}^\flat)\)}
        \label{fig:double_or_nothing}
    \end{subfigure}\\

    \vspace{0.3cm}
    \begin{subfigure}[b]{0.45\textwidth}
        \centering
        \includegraphics[width=0.8\linewidth]{Kamikaze.JPG}
        \caption{Kamikaze}
        \caption*{\(K = (\mathfrak{u}^\sharp, \mathfrak{u}^\flat)\)}
        \label{fig:kamikaze}
    \end{subfigure}
    \qquad
    \begin{subfigure}[b]{0.45\textwidth}
        \centering
        \includegraphics[width=0.5\linewidth]{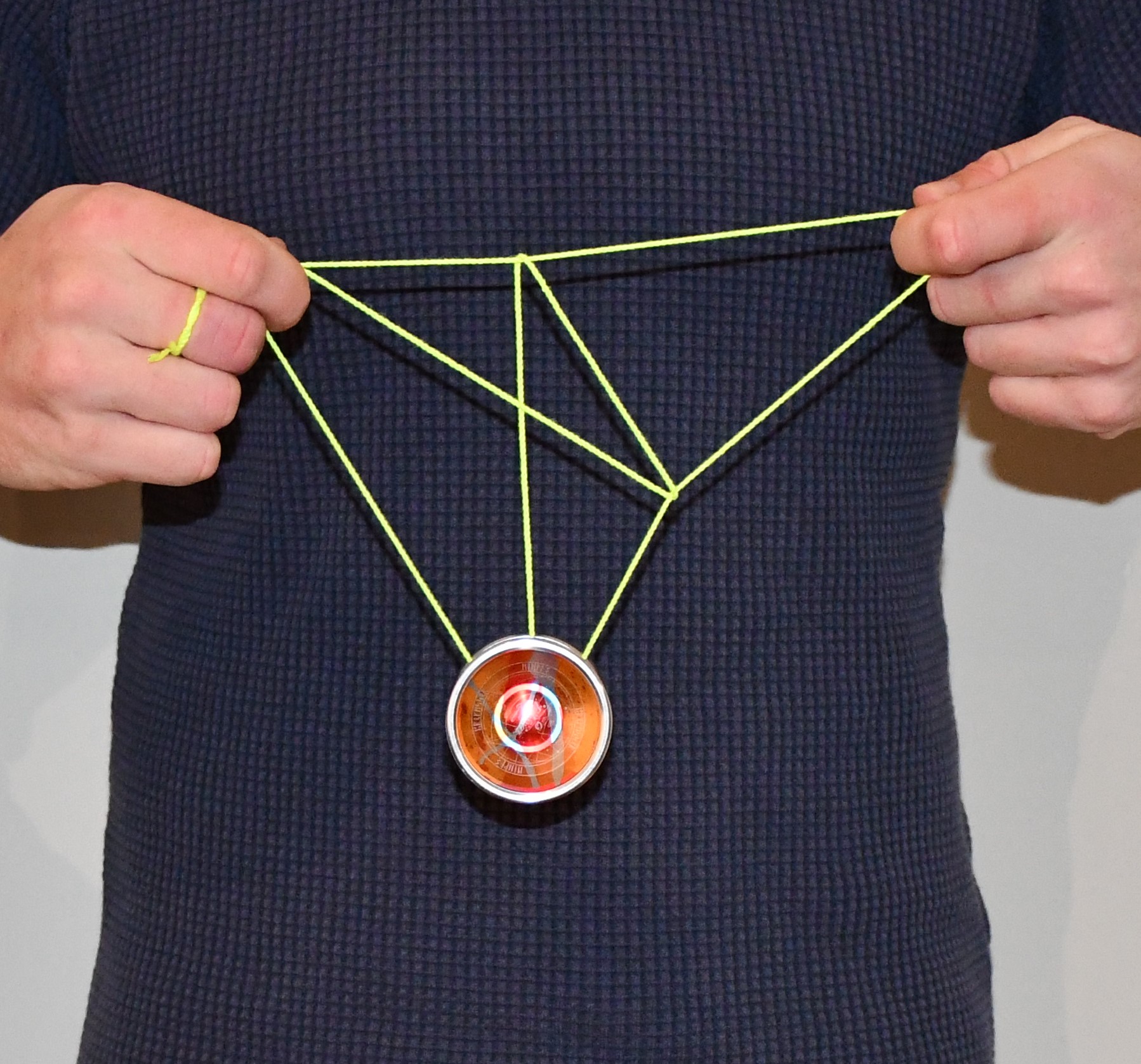}
        \caption{Tower}
        \caption*{\(To = (\mathfrak{u}^\sharp, \mathfrak{u}^\flat)\)}
        \label{fig:tower}
    \end{subfigure}\\

    \vspace{0.3cm}
    \begin{subfigure}[b]{0.45\textwidth}
        \centering
        \includegraphics[width=0.6\linewidth]{Forward_GT.JPG}
        \caption{Green Triangle}
        \caption*{\(GT^\sharp = (\mathfrak{u}^\sharp, 3_1^\flat)\) pictured}
        \caption*{\(GT^\flat = (3_1^\sharp, \mathfrak{u}^\flat)\)}
        \label{fig:green_triangle}
    \end{subfigure}
    \qquad
    \begin{subfigure}[b]{0.45\textwidth}
        \centering
        \includegraphics[width=0.45\linewidth]{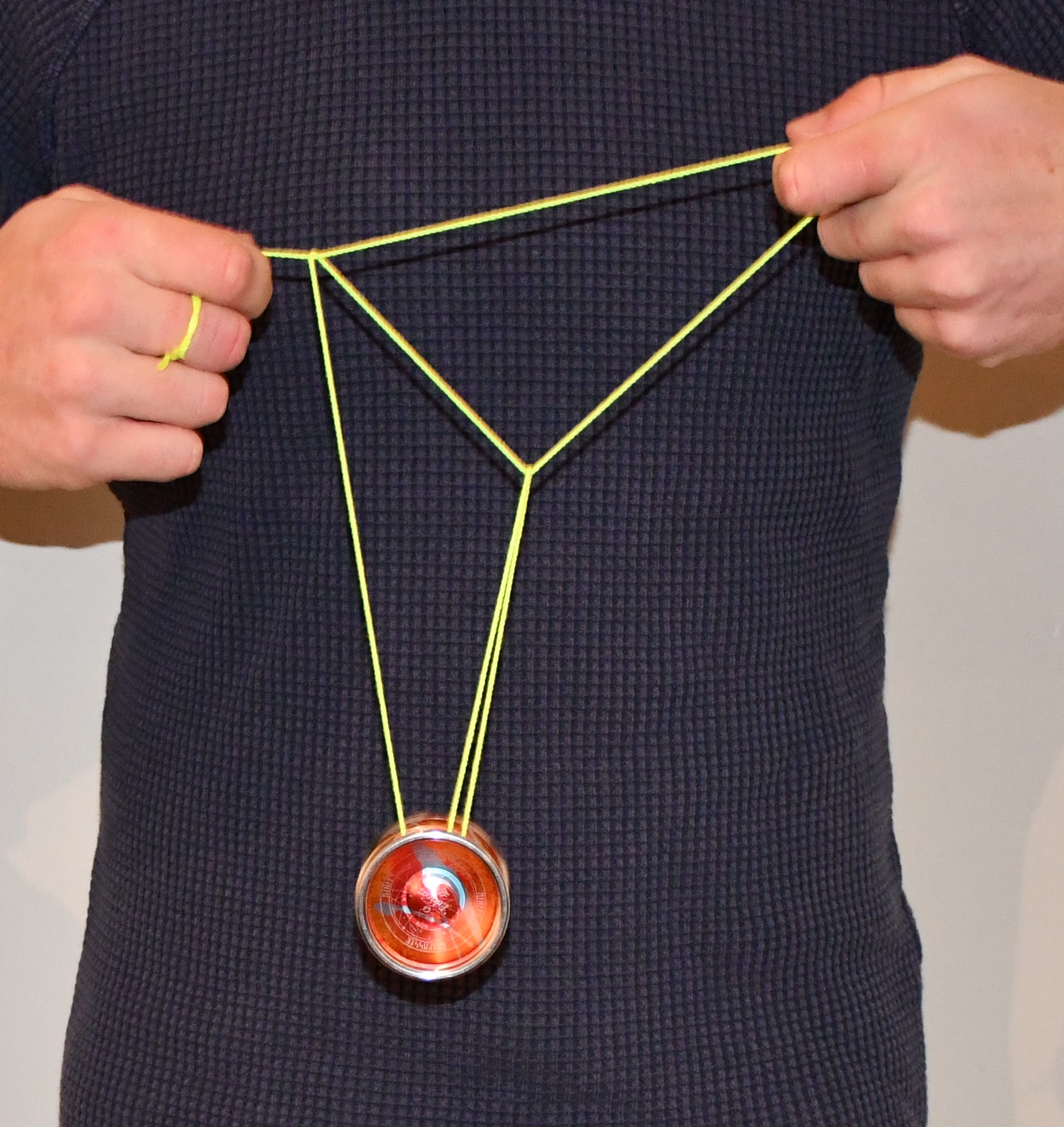}
        \caption{Double Green Triangle}
        \caption*{\(DGT^\flat = (4_1^\sharp, \mathfrak{u}^\flat)\) pictured}
        \caption*{\(DGT^\sharp = (\mathfrak{u}^\sharp, 4_1^\flat)\)}
        \label{fig:double_green_triangle}
    \end{subfigure}\\

    \vspace{0.3cm}
    \begin{subfigure}[b]{0.45\textwidth}
        \centering
        \includegraphics[width=0.45\linewidth]{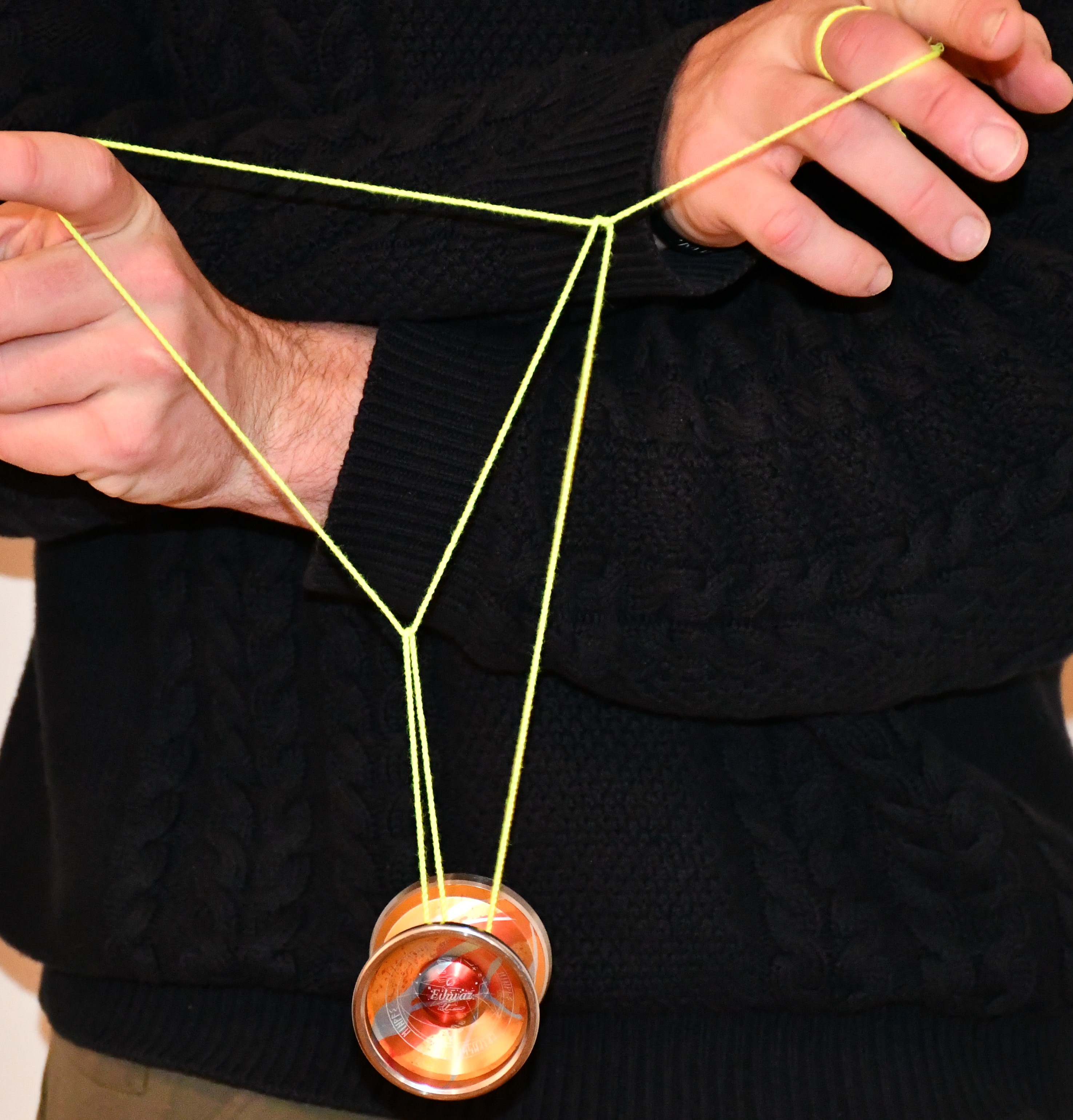}
        \caption{Double Green Triangle Variation}
        \caption*{\(DGT_\xi = (3_1^\sharp, 5_2^\flat)\)}
        \label{fig:double_green_triangle_var}
    \end{subfigure}
    \hfill
    \begin{subfigure}[b]{0.45\textwidth}
        \centering
        \includegraphics[width=0.65\linewidth]{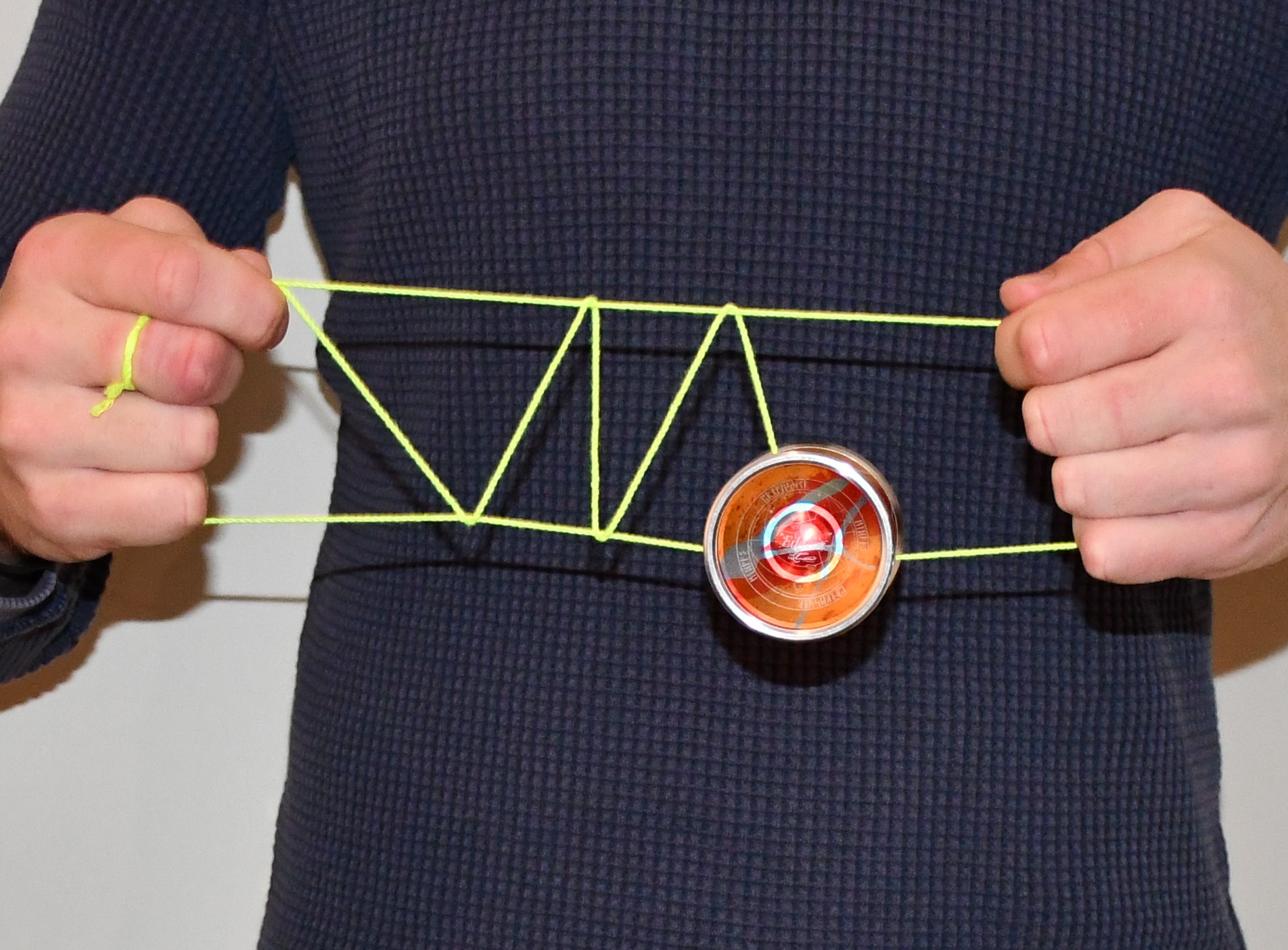}
        \caption{Ladder}
        \caption*{\(L_{4\Delta} = (6_2^\sharp, \mathfrak{u}^\flat)\)}
        \label{fig:ladder}
    \end{subfigure}
    \caption{Common beginner and more intricate mounts, with images of those mounts in play and their corresponding knot-theoretic definitions: (\protect\subref{fig:trapeze}) Trapeze, (\protect\subref{fig:double_or_nothing}) Double or Nothing, (\protect\subref{fig:kamikaze}) Kamikaze, (\protect\subref{fig:tower}) Tower, (\protect\subref{fig:green_triangle}) Green Triangle, (\protect\subref{fig:double_green_triangle}) Double Green Triangle, (\protect\subref{fig:double_green_triangle_var}) Double Green Triangle Variation, and (\protect\subref{fig:ladder}) Ladder. For brevity, where a mount has forward and backward-hopped variants, only one is pictured.}
    \footnotemark
    \label{Mount Table}
\end{figure}

We describe a selection of 1A mounts in terms of prime knots, generally ordered by their number of crossings (except `Magic Knots', which are discussed separately in~\ref{Magic Knots}). Some mounts are analysed in detail, followed by a comprehensive listing of all evaluated mounts in Figure~\ref{Mount Table}.  
Most beginner mounts correspond to the unknot in both forward- and backward-hopped senses. Examples include the `Trapeze', the `Double or Nothing', and the `Kamikaze' mount. The construction of forward-hopped and backward-hopped knots for the Kamikaze mount (\(K = (\mathfrak{u}^\sharp, \mathfrak{u}^\flat)\)) is illustrated in Figure~\ref{Kamikaze}.

\begin{figure}[h!]
    \centering
    \begin{subfigure}{\textwidth}
        \centering
        \includegraphics[width=0.5\textwidth]{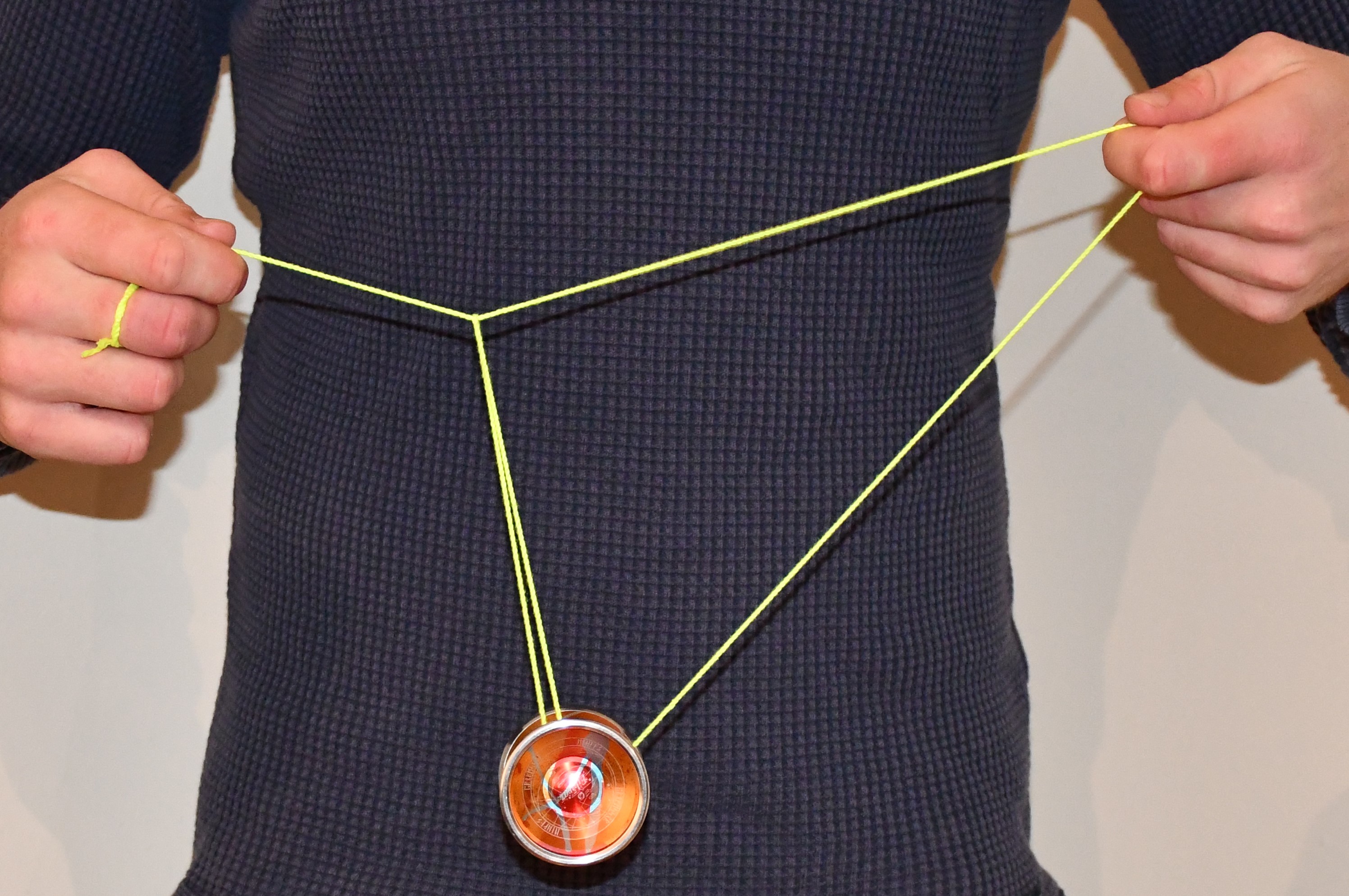}
        \caption{The GT$^\sharp$ mount in play}
        \label{fig:gt_in_play}
    \end{subfigure}\\
    \vspace{.2cm} 
    \begin{subfigure}[b]{0.35\textwidth}
        \centering
        \includegraphics[width=\textwidth]{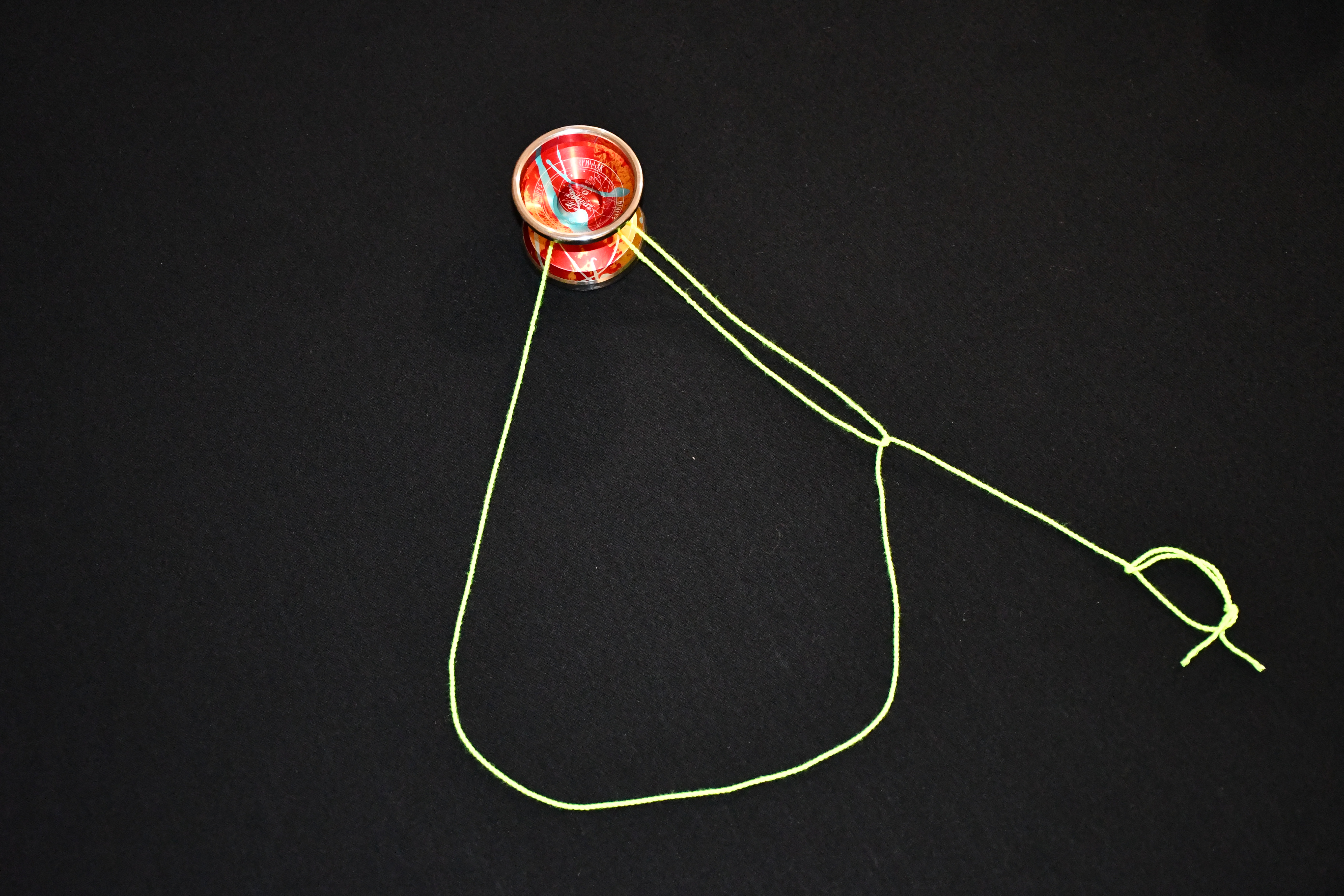}
        \caption{Mount lay flat with a stationary Yo-Yo}
        \label{fig:gt_flat}
    \end{subfigure}
    \quad
    \begin{subfigure}[b]{0.35\textwidth}
        \centering
        \includegraphics[width=\textwidth]{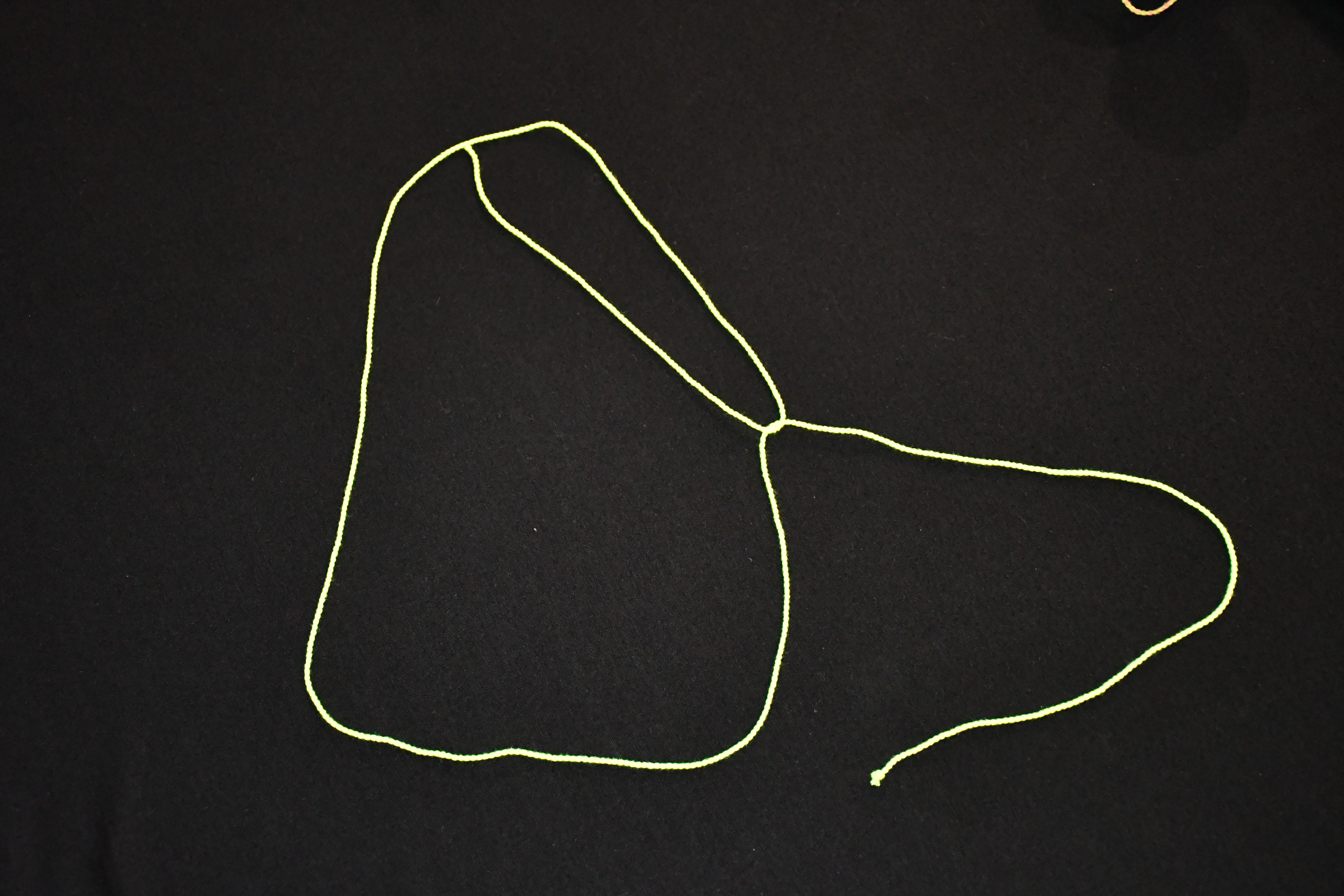}
        \caption{The mount with the Yo-Yo removed}
        \label{fig:gt_no_yoyo}
    \end{subfigure}
    \\
     \vspace{.2cm} 
     \begin{subfigure}[b]{0.35\textwidth}
        \centering
        \includegraphics[width=\textwidth]{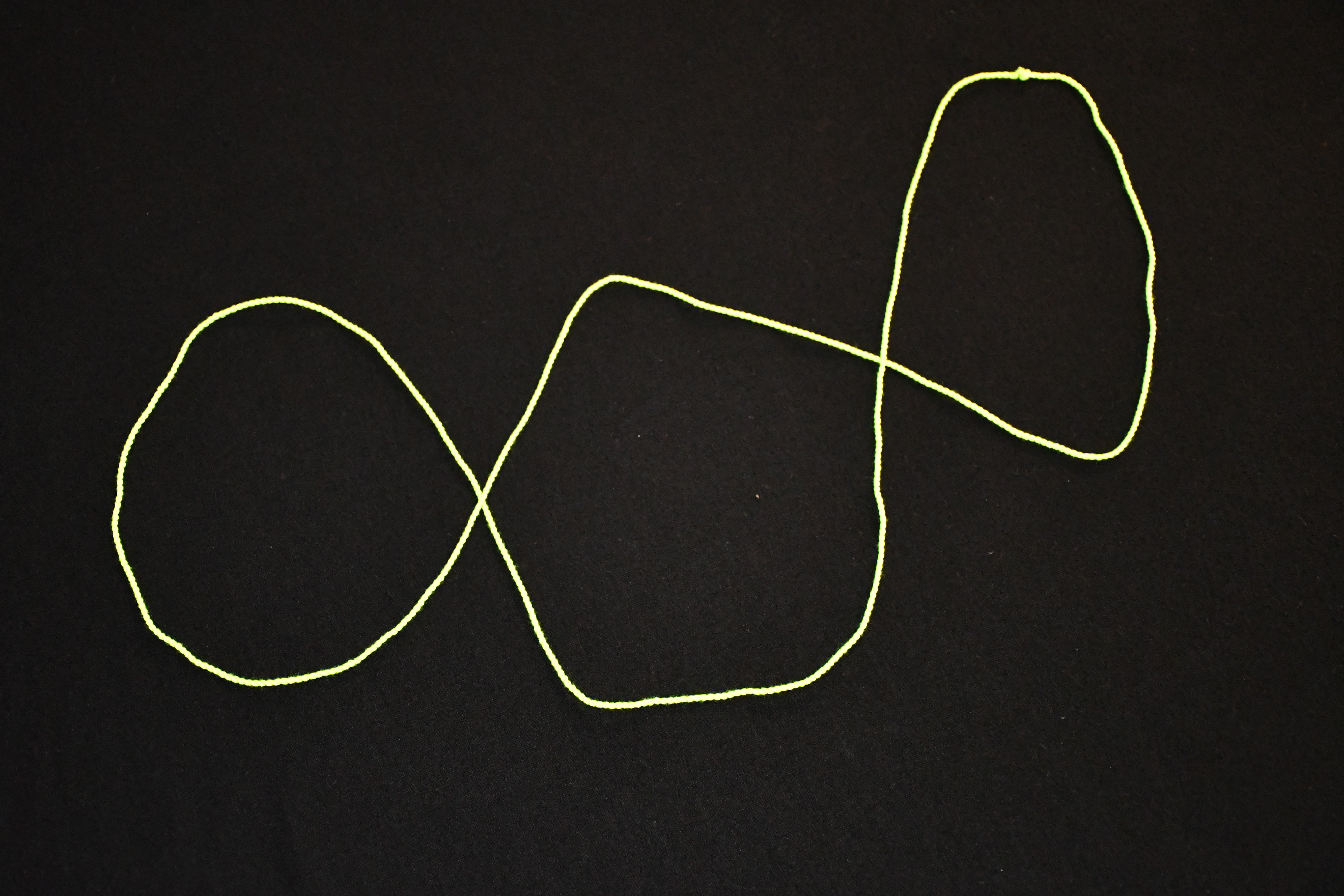}
        \caption{The forward-hopped knot}
        \label{fig:gt_forward}
    \end{subfigure}
\quad    \begin{subfigure}[b]{0.35\textwidth}
        \centering
        \includegraphics[width=\textwidth]{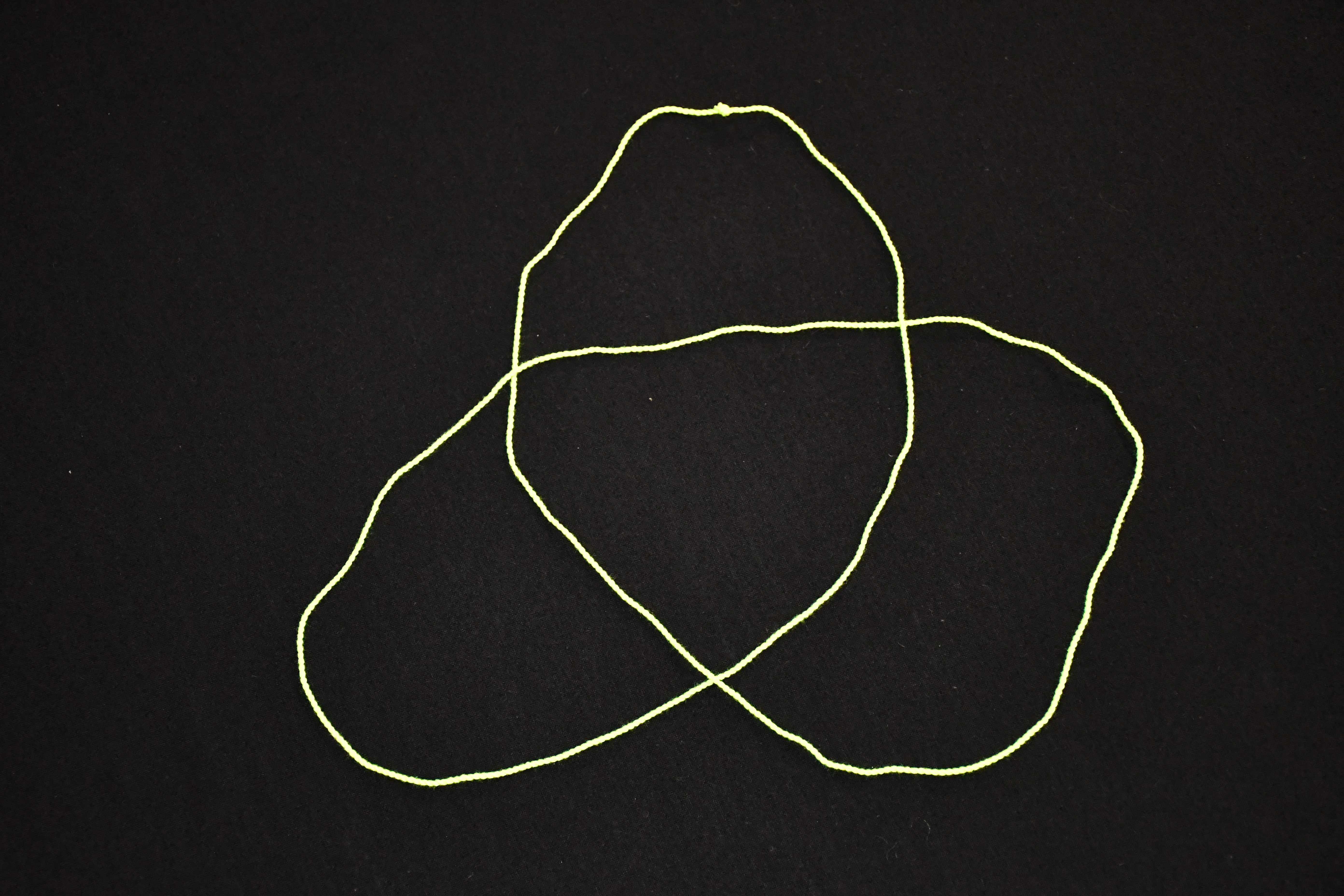}
        \caption{The backward-hopped knot}
        \label{fig:gt_backward}
    \end{subfigure}
    \caption{Construction of the forward- and backward-hopped knots for the Green Triangle (forward-hopped) mount. The top image (\subref{fig:gt_in_play}) shows the mount in play. The bottom row illustrates the key steps in the process: the mount lying flat with a stationary Yo-Yo (\subref{fig:gt_flat}), the mount with the Yo-Yo removed and the slipknot replaced (\subref{fig:gt_no_yoyo}), the resulting forward-hopped knot (\subref{fig:gt_forward}), and the backward-hopped knot (\subref{fig:gt_backward}). The forward-hopped knot recovers the unknot with two Reidemeister moves, while the backward-hopped knot forms the trefoil ($3_1$).}
    \label{GT}
\end{figure}

The most basic mount that is not a double of unknots is known in Yo-Yo parlance as a `Green Triangle' (see Figure~\ref{GT}). A Green Triangle (hereafter GT) refers to either of the mounts 
\((\mathfrak{u}^\sharp, 3_1^\flat)\) or \((3_1^\sharp, \mathfrak{u}^\flat)\). These mounts are obtained differently in play and require distinct methods to `dismount'. Specifically, we define
\begin{equation}
GT^\sharp = (\mathfrak{u}^\sharp, 3_1^\flat) \hspace{0.5cm} \textrm{and} \hspace{0.5cm} GT^\flat = (3_1^\sharp, \mathfrak{u}^\flat)
\end{equation}
as the forward-hopped GT and backward-hopped GT, respectively.  

When the string formation is held parallel to and in front of the torso in side-style play, a GT may be dismounted—meaning the unknot may be recovered—by hopping either `out the front' or `out the back', depending on whether the Yo-Yo is in a \(GT^\sharp\) or \(GT^\flat\) mount. Hopping in the wrong direction leaves the tied knot intact, which is typically undesirable unless transitioning to a more complex mount.  
The construction of the forward-hopped and backward-hopped knots for \(GT^\sharp\) is illustrated in Figure~\ref{GT}. The backward-hopped trefoil obtained from \(GT^\sharp\) is right-handed, while the forward-hopped trefoil from \(GT^\flat\) is left-handed.
\newpage
\newpage
\newpage

\footnotetext{The ladder in Figure~\ref{fig:ladder} forms four triangles when held in play. The ladder can be constructed with (in theory) as many triangles as a player wishes (minimum 2), and this choice will result in a different knot-theoretic mount, but will always result in the unknot when backward-hopped.}

\newpage

One important feature to highlight from these results is that whilst it is clear that some mounts can appear very different and yet be equivalent, it is also true that two mounts which look identical unless closely inspected can be topologically distinct in important ways. For example, \(DGT^\flat\) and \(DGT_\xi\) are very difficult to tell apart from a distance (except for being mirror images of one another), and yet admit entirely distinct mount definitions. This is because the string segment hanging off of the top string in \(DGT_\xi\) is wrapped twice, whereas it is only wrapped once in \(DGT^\flat\), meaning the unknot cannot be recovered from \(DGT_\xi\) with a single hop and must be dismounted in multiple stages.

\subsubsection{Magic Knots}
\label{Magic Knots}

It may seem surprising that the `Tower' mount displayed in Figure~\ref{Mount Table} is the unknot in both the forward and backward sense. This is an example of a `magic knot'. Since there are infinitely many knots equivalent to \(\mathfrak{u}\)\footnote{For example, one can continuously apply Reidemeister move I to the string to obtain a new diagram each time the move is applied.}, it is possible to construct elaborate mounts that are nonetheless \((\mathfrak{u}^\sharp, \mathfrak{u}^\flat)\). The creation of such mounts and the subsequent revelation that the string arrangement does not result in a nontrivial knot underpins the `magic knot' effect. 

\begin{figure}[h!]
\centering
\begin{minipage}{0.45\textwidth}
\includegraphics[width=0.9\textwidth]{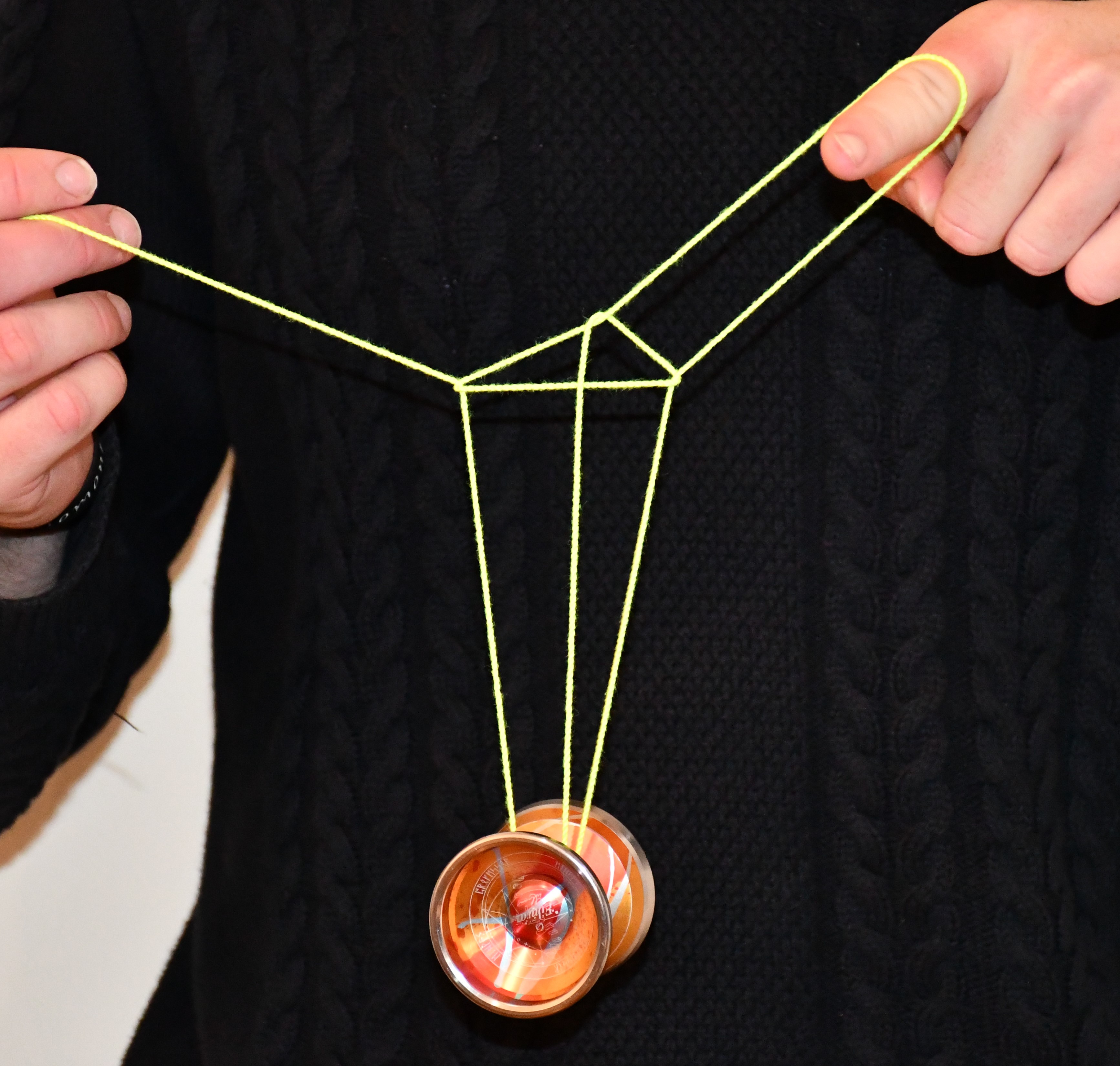}
\caption{The in-play \(GT^\sharp\circ To\) mount.}
\label{GTTower}
\end{minipage}\hfill
\begin{minipage}{0.45\textwidth}
\includegraphics[width=0.9\textwidth]{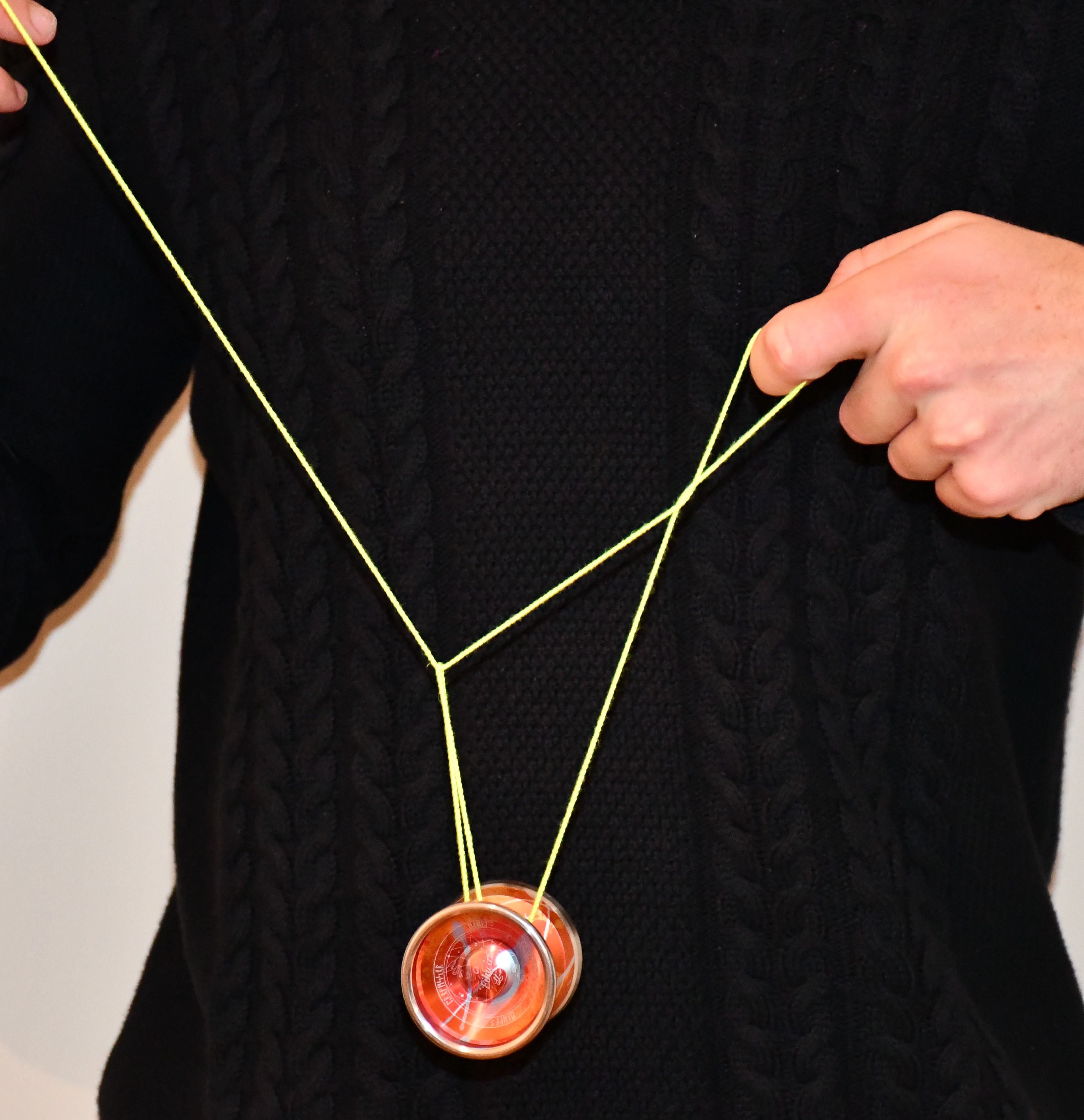}
\caption{The resultant mount from the `Brent Stole' trick.}
\label{BrentStole}
\end{minipage}

\end{figure}

New magic knot like mounts can be built using the techniques established in this paper. For instance, since \(To\) is the unknot in both the forward and backward sense, inserting this arrangement into any component of a mount does not change its crossing number. Consequently, a mount constructed by inserting a \(To\) mount into a \(GT^\sharp\) mount will have the same knot-theoretic definition as \(GT^\sharp\) itself. Such a mount, pictured in Figure~\ref{GTTower}, recovers the unknot with a single hop out the front, despite its seemingly complex appearance. 

\begin{definition}[Mount Composition]
We compose mount \(M_1\circ M_2\) by traversing from the origin to the terminus of the string, first adding all oriented crossings from the in-play mount \(M_1\), and then adding all crossings from the in-play mount \(M_2\), with the Yo-Yo resting on the string segment as it would in \(M_2\).
\end{definition}

While analogous to the knot sum, mount composition does not generally satisfy
\begin{equation}
M_1\circ M_2 \ne (M_1^\sharp \# M_2^\sharp, M_1^\flat \# M_2^\flat),
\end{equation}
because in-play mounts are `undecided' regarding which knot they will form. For example, a \(GT^\sharp\) mount shows two crossings in play—one over and one under—from the player's viewpoint, which is insufficient to form a trefoil until hopping backward introduces the necessary crossing. 

As an illustration, for two identical \(GT^\sharp\) mounts:
\begin{equation}
GT^\sharp \circ GT^\sharp = (3_1^\sharp, [3_1 \# 3_1]^\flat).
\end{equation}
Mirrored mounts can `cancel' each other, as in
\begin{equation}
GT^\sharp \circ GT^\flat = (\mathfrak{u}^\sharp, \mathfrak{u}^\flat).
\end{equation}

We leave a rigorous study of mount composition for future work. A key property useful in constructing magic knots is:
\begin{equation}
M \circ (\mathfrak{u}^\sharp, \mathfrak{u}^\flat) = M, \quad (\mathfrak{u}^\sharp, \mathfrak{u}^\flat) \circ M = M.
\end{equation}

\subsection{Trick Elements Describable by Reidemeister Moves}

It is instructive to introduce this discussion with an example. Consider the `Double or Nothing' (\(D = (\mathfrak{u}^\sharp, \mathfrak{u}^\flat)\)) and the `Tower' (\(To = (\mathfrak{u}^\sharp, \mathfrak{u}^\flat)\)) mounts. These two mounts are equivalent in the knot-theoretic sense but differ in visual presentation, as shown in Figure~\ref{Mount Table}. One can obtain the Tower from the Double or Nothing in a single move—by hopping the Yo-Yo off the front string segment over the middle segment and landing on the segment closest to the player's torso.  

If we form the backward-hopped knot for both mounts, the results differ by a single Reidemeister move. Specifically, this trick element is describable by Reidemeister move II (moving one loop completely over another), which does not alter the writhe of the mount. The only move that affects the writhe is Reidemeister move I (twisting), which has important consequences for string torsion\footnote{String torsion is universally mislabelled as `string tension' in the Yo-Yo community. We adopt the nomenclature `string torsion' to clarify the actual effect of twists on the string.}.  

We define the writhe of a mount as the writhe of the diagram representing the mount laid flat with the Yo-Yo removed, as in Figure~\ref{fig:gt_no_yoyo}. This provides the most accurate representation of how writhe impacts in-play mechanics. For example, a mount nearly identical to \(GT^\sharp\) arises from a trick known as the `Brent Stole'. This trick produces a GT with two additional twists from a `sleeping' (fully un-mounted) position. The resultant mount, displayed in Figure~\ref{BrentStole}, changes the writhe as follows:
\begin{equation}
\textrm{Wr}(GT^\sharp) = 0, \hspace{0.5cm} \textrm{Wr}(BS) = -2,
\end{equation}
where \(BS = (\mathfrak{u}^\sharp, 3_1^\flat)\) denotes the mount resulting from the Brent Stole trick. This difference affects subsequent play unless the twists are undone, corresponding to two instances of Reidemeister move I.

By identifying trick elements representable by Reidemeister moves, players can anticipate whether subsequent components of a combination require additional care to dismount. Furthermore, by analysing the writhe of a mount (induced by Reidemeister move I elements), one can flag significant string torsion accumulated during a trick and compensate accordingly.  

\section{Conclusion and Outlook}

This paper lays the groundwork for using modern Yo-Yo play to provide physical realisations of knots, via appropriate post-processing procedures. Simultaneously, it establishes a systematic framework to categorise mounts in practical Yo-Yo play, enabling players to construct new tricks with confidence regarding how they will be dismounted or progressed.  

With tangible links to embeddings, orientations, and knot equivalence, the constructions presented here also offer a useful pedagogical tool for introducing knot theory and topology in educational contexts. We emphasise that this work represents a first step. A natural extension in the 1A context is to examine surgeries in sophisticated tricks; for instance, a Yo-Yo may be popped into the air and the string whipped around it to instantly form a perfect Double Green Triangle.  

Players also exploit physics in advanced play. In the `tension hook' trick, the string is held taut and then released, whipping the string to recapture the Yo-Yo mid-air in a desired mount. These examples suggest avenues to explore the intersection of continuum mechanics and knot theory, examining the behaviour of a one-dimensional embedding (possibly treated via fibre mechanics) in ambient 3-space according to the manifold's mechanical properties.  

\begin{figure}[h!]
\caption{The in-play `Kink Trapeze' 3A mount}
\centering\includegraphics[scale=0.18]{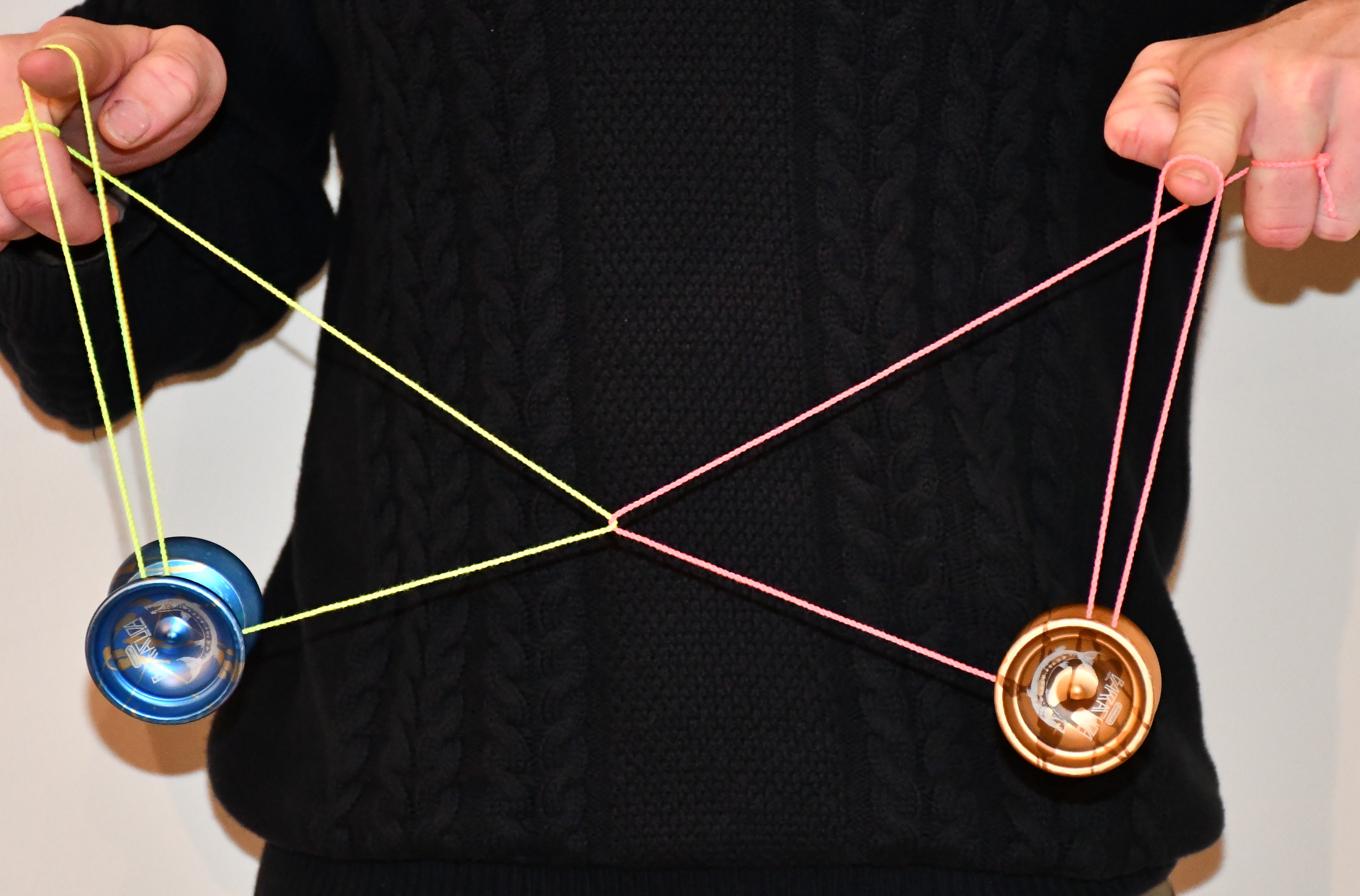}
\label{3A}
\end{figure}

The theory can also be expanded to encompass links and braids, particularly in the context of 3A play. For example, by modifying the double mount shown in Figure~\ref{3A} and applying the mount modification procedure described in Section~\ref{Setup} to \textit{both} Yo-Yos, one obtains the Hopf link~\cite{KnotBook}. Such extensions highlight the potential for further research connecting competitive Yo-Yo play, topology, and applied mechanics.  

Thus, it is hoped that this work will initiate a dialogue between the Yo-Yo and knot theory communities, opening avenues for more sophisticated and surprising tricks in the ever-evolving frontier of modern Yo-Yo play, while also providing inspiration for, and physical realisations of, state-of-the-art results in the mathematics of knots and links.
\label{Conclusion}

\section*{Declarations: Funding and/or Conflicts of interests/Competing interests}

The authors declare no conflicts of interest.\\

The Yo-Yos used for all 1A and 3A mounts displayed in this work were the C3YoyoDesign Eihwaz and the Duncan Barracuda, respectively. The string used was C3YoyoDesign Pro-String.\\

No AI tools were used in this work.


\footnotesize
    
\bibliographystyle{unsrt}
\bibliography{bib}

\end{document}